 \documentclass[final,1p,times]{elsarticle}

\usepackage{amssymb}
\usepackage{amsthm}
\usepackage[centertags]{amsmath}
\usepackage{pinlabel-test} 

\theoremstyle{plain}
\begingroup
\newtheorem{thm}{Theorem}[section]

\endgroup

\theoremstyle{remark}
\begingroup
\newtheorem{rem}[thm]{Remark}
\endgroup

\theoremstyle{definition}
\begingroup

\endgroup

\begin{document}

\begin{frontmatter}

\title{Mechanics of dislocation pile-ups: \\A unification of scaling regimes}

\author{L. Scardia\corref{cor1}\fnref{label1,label2,label3}}
\ead{L.Scardia@tue.nl}
\author[label2]{R.H.J. Peerlings}
\ead{R.H.J.Peerlings@tue.nl}
\author[label3,label4]{M.A. Peletier}
\ead{M.A.Peletier@tue.nl}
\author[label2]{M.G.D. Geers}
\ead{M.G.D.Geers@tue.nl}

\address[label1]{Materials innovation institute (M2i)}
\address[label2]{Department of Mechanical Engineering, Technische Universiteit Eindhoven, The Netherlands}
\address[label3]{Department of Mathematics and Computer Sciences, Technische Universiteit Eindhoven, The Netherlands}
\address[label4]{Institute for Complex Molecular Systems, Technische Universiteit Eindhoven, The Netherlands}

\cortext[cor1]{Corresponding author. Tel.: +31402472644; Fax: +31402447355}

\begin{abstract}
This paper unravels the problem of an idealised pile-up of $n$ infinite, equi-spaced walls of edge dislocations at equilibrium.
We define a dimensionless parameter that depends on the geometric, constitutive and loading parameters of the problem, and we
identify five different scaling regimes corresponding to different values of that parameter for large $n$.
For each of the cases we perform a rigorous micro-to-meso upscaling, and we obtain five expressions for the mesoscopic (continuum) internal stress.
We recover some expressions for the internal stress that are already in use in the mechanical community, as well as some new models.
The results in this paper offer a unifying approach to such models, since they can be viewed as the outcome of the same
discrete dislocation setup, for different values of the dimensionless parameter (i.e., for different local dislocations arrangements).
In addition, the rigorous nature of the upscaling removes the need for ad hoc assumptions.
\end{abstract}


\end{frontmatter}

\section{Introduction}
\label{Intro}
Dislocations occupy a central position in discussions of the permanent deformation of metals because of their role as the main carriers of plastic
deformation. Therefore it is necessary to incorporate their presence, or the main effect of their presence, in a plasticity model that aims for
a predictive power. However, since the typical number of dislocations even in a small sample of the material is very high, formulating a model that keeps track of every
single dislocation is out of reach except for very small-scale problems. This explains the interest in describing the collective behaviour of
dislocations in terms of a continuum quantity: the dislocation density.

The challenge in this scale transition consists in describing the time evolution of the dislocation density in a physically-driven way, by performing a
rigorous upscaling from the dislocation scale to the dislocation density scale (also called \textit{meso-scale}).
This task has been mainly pursued phenomenologically or by means of a statistical mechanics approach and has produced a number of competing models
(e.g.~\cite{DengEl-Azab09,El-Azab00,EvGeers,Groma97,GromaBalogh99,Groma,LimkumnerdVan-der-Giessen08,RoyAcharya06}). In the case of parallel, straight dislocations, the evolution of the dislocation density $\rho(x,t)$ (here we assume, for simplicity, that all the dislocations have the same Burgers vector) is described in terms of a continuity equation of the following form:
$$
\partial_ t\rho + \partial_x(\rho v) = 0,
$$
where the velocity $v$ driving the evolution is of the form
\begin{equation}\label{form:v}
v = \frac{b}{B}\,(\sigma-\sigma_{\textrm{\scriptsize{int}}}),
\end{equation}
with $b$ denoting the Burgers vector, $B$ a linear drag coefficient, $\sigma$ an externally applied shear stress (constant, for simplicity) and $\sigma_{\textrm{\scriptsize{int}}} = \sigma_{\textrm{\scriptsize{int}}}(x,\rho,\partial_x \rho)$ an internal stress accounting for the net effect of the interactions among dislocations.
The different models available in the engineering literature differ in the expressions of the internal stress that they propose; the range of validity of the proposed theories is typically unclear, as well as the conditions under which one of them is to be preferred over another one. In fact most of the models have been derived phenomenologically, often starting from an \emph{ad hoc} Ansatz, or assumptions have been made in the derivation that are not always justified.

\medskip

Here comes the essential difference of our approach: we obtain a continuum model for the dislocation density from a more fundamental discrete dislocation model using a rigorous mathematical approach ($\Gamma$-convergence). As a consequence, we obtain an exact classification of limiting behaviour for the system that we study, which unifies existing, independently derived descriptions into a coherent framework, and identifies new regimes that have not yet been studied.

In order to make each step in the derivation completely justified and rigorous, our starting point has to be an idealised dislocation configuration. More precisely, we consider the discrete model of an idealised pile-up of dislocations studied in \cite{RPG} (see also \cite{TRB,Hall10, Hall11,BaskMes10, MesarovicBaskaranetc10}). This model describes the equilibrium positions of $n$ dislocation walls under the influence of an applied stress $\sigma$ that pushes the walls towards an impenetrable barrier; the barrier is modelled as a wall of pinned dislocations at $x_0=0$ (see Section \ref{Prob-Statement} for the detailed description).

The discrete equilibrium equations for the positions of the $n$ walls can be written in the general form
\begin{equation}\label{discretesigma}
\sigma_{\textrm{\scriptsize{int}}}^i - \sigma=0, \quad i=1,\dots,n,
\end{equation}
where the discrete internal stress for the $i$-th wall is the sum of the contributions due to the interactions with the other walls (see (\ref{Eq:1}) for the detailed expression of the equations). Therefore the system (\ref{discretesigma}) is nothing but $v_i=0$, where $v_i$ is the velocity of the $i$-th wall.
At this point one can intuitively imagine that passing to the (continuum) limit in the discrete equation (\ref{discretesigma}) should give a continuum analog of (\ref{form:v}) for $v=0$, and consequently an expression for a continuum version of the internal stress. This is exactly the object we want to characterise.

Our approach is different from the upscaling procedure followed in the quoted papers that considered the same pile-up configuration. In \cite{Hall10, Hall11} the convergence of stationary states is proved using formal methods, in a special case of our analysis. Moreover, while \cite{BaskMes10, MesarovicBaskaranetc10} perform a two-step discrete-to-continuum upscaling by smearing out the dislocations first in the slip plane and then in the in-wall direction, we upscale in the two directions simultaneously. We stress that the continuum model we derive is in perfect agreement with our starting discrete model. We recall that in \cite{BaskMes10, MesarovicBaskaranetc10}, as the authors point out, the resulting continuum model has to be corrected to incorporate the missing interaction.

\medskip

\begin{figure}[htb]
\begin{center}$
\begin{array}{ll}
\includegraphics[width=1.8in]{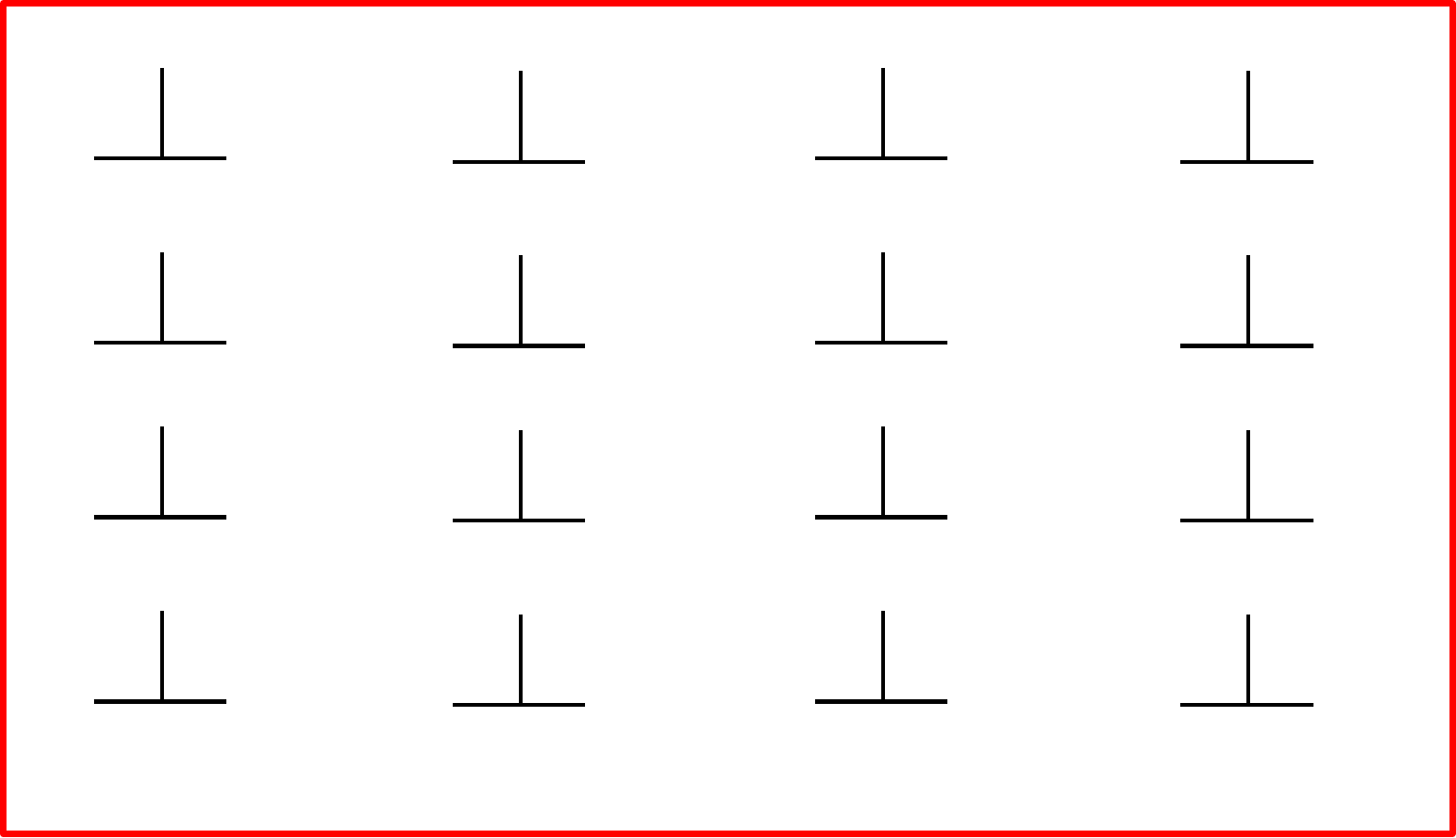} & \hspace{2.3cm}
\includegraphics[width=1.8in]{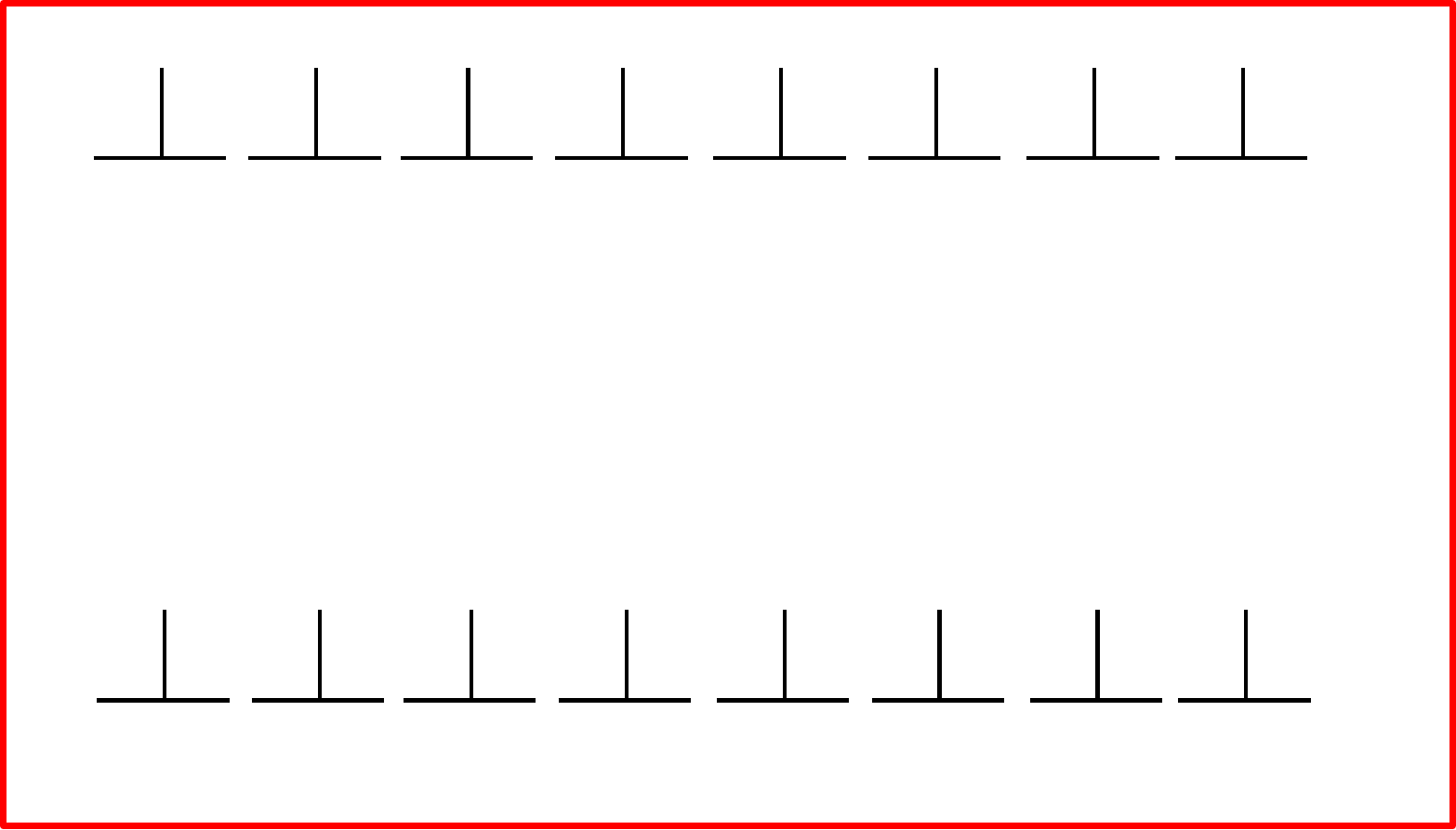}
\end{array}$
\end{center}
\caption{Two different arrangements corresponding to the same density in the rectangle.}
\label{Densitypattern}
\end{figure}

An interesting novelty in our result is that the mesoscopic internal stress we obtain does not depend on the dislocation density (and its gradient) only, but it also contains some more \textit{local} information about the discrete arrangement of the dislocations, that the density alone would fail to capture.
To be more precise, consider the two arrangements in Figure \ref{Densitypattern}. They correspond to the same density, i.e., to the same number of dislocations per unit volume. But the dislocation patterns in the two cases are different and this will result in different expressions for the upscaled internal stress describing the overall interactions in the two cases.

This additional \textit{local} information on the arrangement of the dislocations can be expressed in terms of an aspect ratio $a$, defined as the ratio between the typical distance $\Delta x$ between two consecutive dislocations in the same row (slip plane) and the distance $h$ between two consecutive dislocations in the same wall. In Figure \ref{Densitypattern} the first arrangement corresponds to $a\sim 1$, the second one to $a\ll 1$.

\medskip

The upscaling procedure we adopt for the micro-to-meso upscaling is based on a variational convergence called $\Gamma$-convergence, which is well-known in the mathematical community and has been successfully applied to a variety of problems in materials science, from fracture mechanics to homogenization, from magnetomechanics to dimension reduction. Our approach focusses on the discrete energy of the system of dislocations (minimisers of the energy are exactly the solutions to (\ref{discretesigma})). The discrete-to-continuum upscaling is done by letting the number of dislocations $n$ become infinitely large. According to the different asymptotic behaviour of the aspect ratio $a$ (i.e., according to the \textit{local} distribution of the dislocations), five different expressions for the continuum energy can be derived (we refer to \cite{GPPS_Math} for the details of the mathematical procedure); and, accordingly, five different expressions for the upscaled internal stress.

The results we obtain show that the simplified discrete model taken as a starting point is not \textit{too simple}.
In fact the internal stresses resulting from our derivation are more general than many well-known models proposed in the engineering literature
(see \cite{EFN}, \cite{EvGeers} and \cite{Groma}). More precisely, our result contains as a special case the internal stresses proposed in the quoted papers,
explaining their range of validity. The comparison with previous models will be the subject of Section \ref{Comparison}.

\medskip

The key advantage of our rigorous approach to  upscaling compared to other methods is that it is exact. This means that once a discrete model is chosen (with its simplifications and limitations) the corresponding upscaled continuum model obtained by following our method is uniquely determined.

Moreover, the fact that no simplifying assumptions are made during the derivation (but only at the level of the discrete model) explains why such a rich variety of continuum models can be obtained, even from a relatively simple model. This is not the case for other upscaling methods, where the starting point can be more general, but a number of restrictions and simplifications are typically made in the derivation, affecting the generality of the upscaled models and their range of validity.

\bigskip
The plan of the paper is as follows.
In Section \ref{Prob-Statement} we state the discrete problem and we comment on our special choice for the arrangement of the dislocations. An overview of the methods and of the results is given in Section \ref{Main}, while the details of the upscaling are considered in Sections \ref{Subcritical}--\ref{Supercritical}. A comparison between our results and other well-known models in the engineering literature is the subject of Section \ref{Comparison}, and a number of interesting comments on the the results of this paper are presented in Section \ref{Comments}. Section \ref{Conclusion} contains a short conclusion.

\section{Statement of the problem}
\label{Prob-Statement}

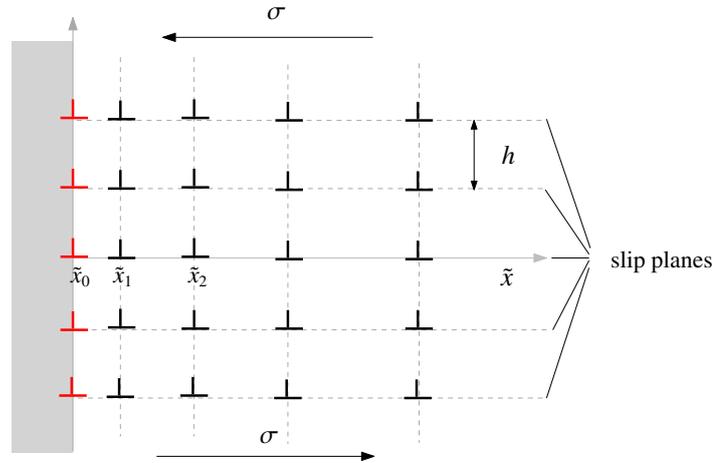
\begin{figure}[htbp]
\labellist
\pinlabel \small$\tilde{x}_0$ at 90 208
\pinlabel \small$\tilde{x}_1$ at 124 208
\pinlabel \small$\tilde{x}_2$ at 184 208
\pinlabel $\tilde{x}$ at 431 208
\pinlabel $h$ [l] at 421 307
\pinlabel $\sigma$ [b] at 247 410
\pinlabel $\sigma$ [t] at 241 90
\pinlabel {\small slip planes} [l] at 508 220
\endlabellist
\begin{center}
\includegraphics[width=3in]{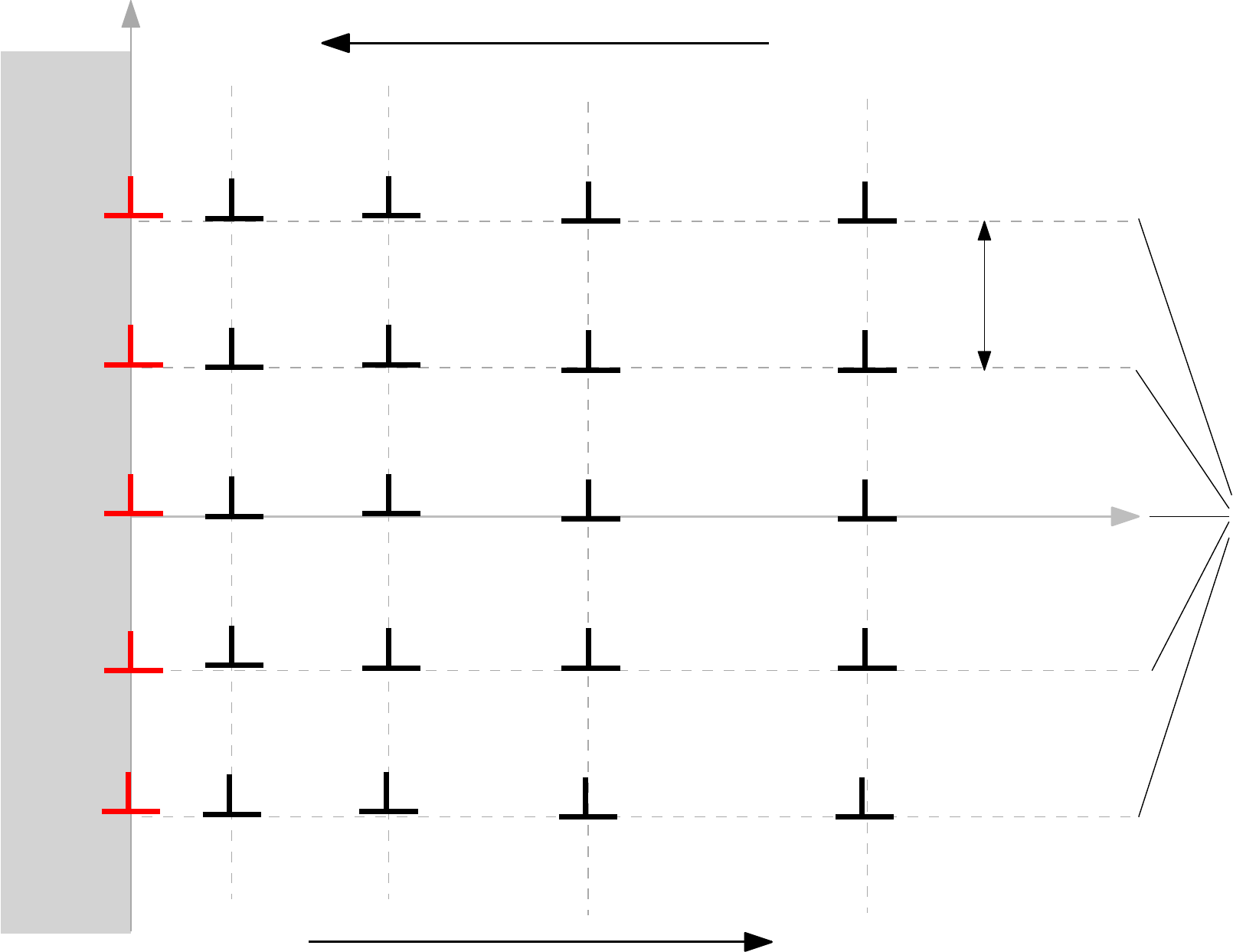}
\end{center}
\caption{Pile-up of discrete dislocation walls. The walls are at positions $\tilde{x}_i>0$ and the barrier is modelled as a wall of pinned dislocations at $\tilde{x} = \tilde{x}_0 =0$.}
\label{Wall}
\end{figure}

The problem we are considering is the equilibrium of a system of $n$ walls of straight edge dislocations (all with the same Burgers vector) under the action of an externally applied shear stress that pushes the walls towards an impenetrable barrier. Consecutive dislocations within the same wall are assumed to be equidistant, at distance $h>0$, and we assume there are infinitely many dislocations in each wall. The impenetrable barrier is modeled as an infinite wall of pinned dislocations (the plane $\tilde{x}=0$). We moreover assume that the walls are perfectly aligned (see Figure \ref{Wall}). Therefore the model is essentially one-dimensional, the only unknowns being the positions $\tilde{x}_1, \dots, \tilde{x}_n$ of the~$n$ dislocation walls.\footnote{From here on, tildes distinguish dimensional quantities from their non-dimensional counterparts; we will define non-dimensional positions $x_i$ below.}

While it is  common in the literature to reduce to the case of straight dislocations and to consider a single slip system (rather than multiple slip or cross-slip), the assumption of having dislocations arranged in vertical walls may sound unnatural (even if used in many papers, see e.g. \cite{BaskMes10, MesarovicBaskaranetc10,Hall10, Hall11,RPG,TRB}). It is clearly a simplification and we are aware of the fact that a complete analysis of a fully 3d-dislocation model would be preferable. On the other hand, this idealised configuration allows us to carry out a rigorous analysis, whereas more random configurations would require a statistical approach. Such a statistical approach, however, would typically require some additional assumptions in the derivation, and therefore it might  well result in a limit model that is less general than the one rigorously derived from a more idealised setting.

Moreover, although the discrete model we consider is highly idealised, it has a number of properties that make it both interesting and not unrealistic. The fact that multiple dislocations move along exactly the same slip plane is natural, because of the way they are generated from \emph{Frank-Read sources} (e.g.~\cite[Sec.~8.6]{HullBacon01}). Also, although the assumption of an arrangement in equispaced vertical walls is clearly an idealisation, it is on the other hand not unrealistic since equispaced vertical walls are minimal-energy configurations. Walls of edge dislocations are locally stable, in the sense that if one of the dislocations deviates from its wall position, either horizontally or vertically, it experiences a restoring force from the other dislocations that pushes it back. Finally, the vertical organization in walls is also justified by correlation functions calculated from numerical simulations~(e.g.~\cite{Groma}).

For these reasons we believe that our simplified discrete system is not too simple and that it is worth exploring it. Moreover, the mesoscopic models we obtain are general enough to contain as a special case several well-known models in the engineering literature, including some models derived by statistical arguments.

\smallskip

The equilibrium positions of the dislocation walls are obtained by solving the following equations for every wall $i=1,\dots,n$:

\begin{equation}\label{Eq:1}
\frac{\pi G b}{2h(1-\nu)}\sum_{\stackrel{j=0}{j\neq i}}^{n}
\varphi \left(\frac{\tilde{x}_i - \tilde{x}_j}{h}\right) - \sigma = 0,
\end{equation}
where $\varphi(s) = \frac{s}{\sinh ^2\pi s}$ is the (globally repulsive) stress governing the wall-wall interactions, $-\sigma$ is the constant applied shear stress, $b$ is the length of the Burgers vector, $\nu$ is the Poisson's ratio and $G$ is the shear modulus.
Setting
\[
K:= \frac{\pi G b}{2(1-\nu)},
\]
we can rewrite (\ref{Eq:1}) as
\begin{equation}\label{Eq:1r}
\frac{K}h\sum_{\stackrel{j=0}{j\neq i}}^{n}
\varphi \left(\frac{\tilde{x}_i - \tilde{x}_j}{h}\right) - \sigma = 0, \quad i=1,\dots,n.
\end{equation}

\medskip

\begin{figure}[htb]
\labellist
\pinlabel {\small (a)} at 180 -15
\pinlabel {\small (b)} at 570 -15
\pinlabel {\small (c)} at 950 -15
\endlabellist
\begin{center}$
\begin{array}{lll}
\includegraphics[width=1.6in]{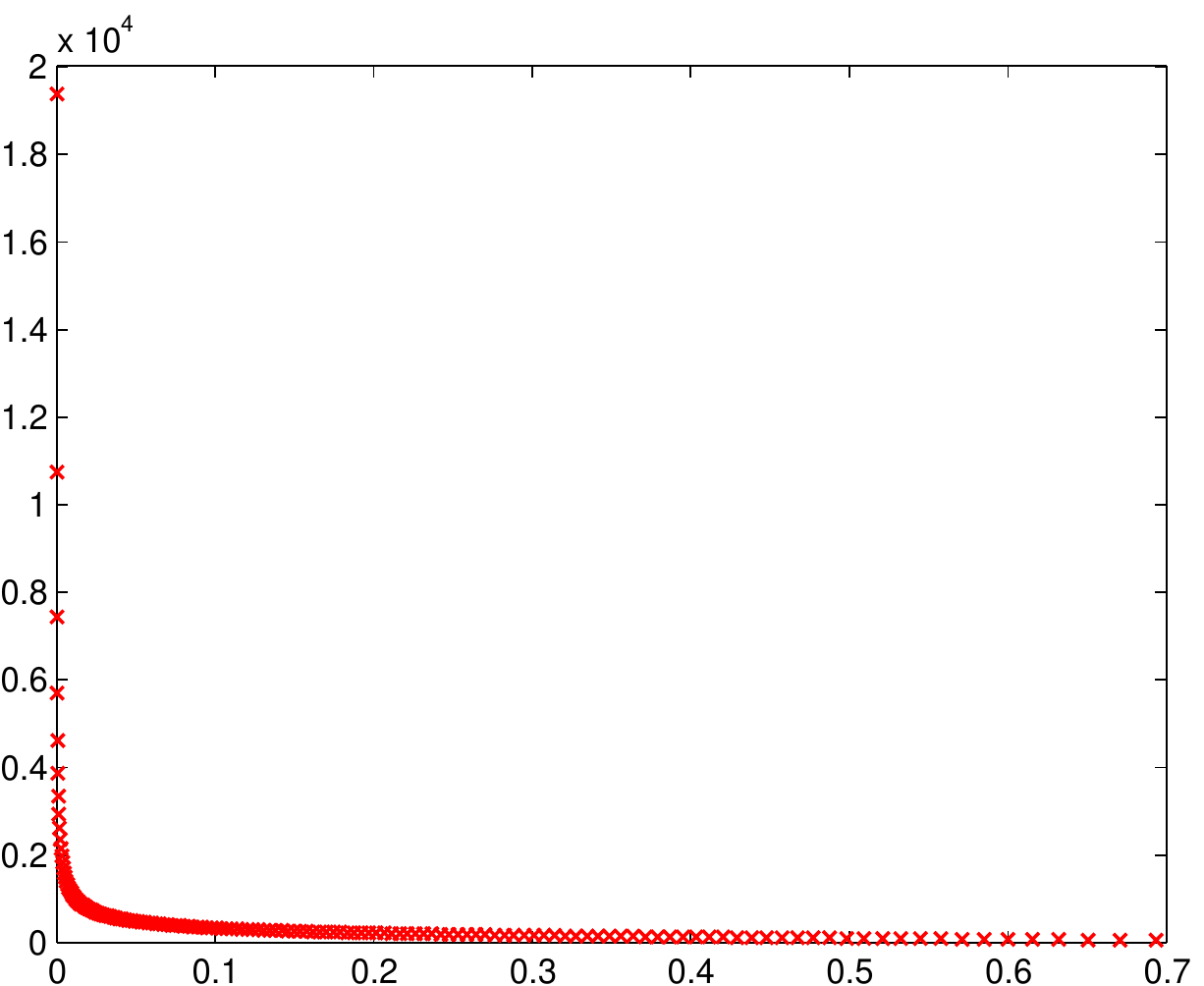} &
\includegraphics[width=1.6in]{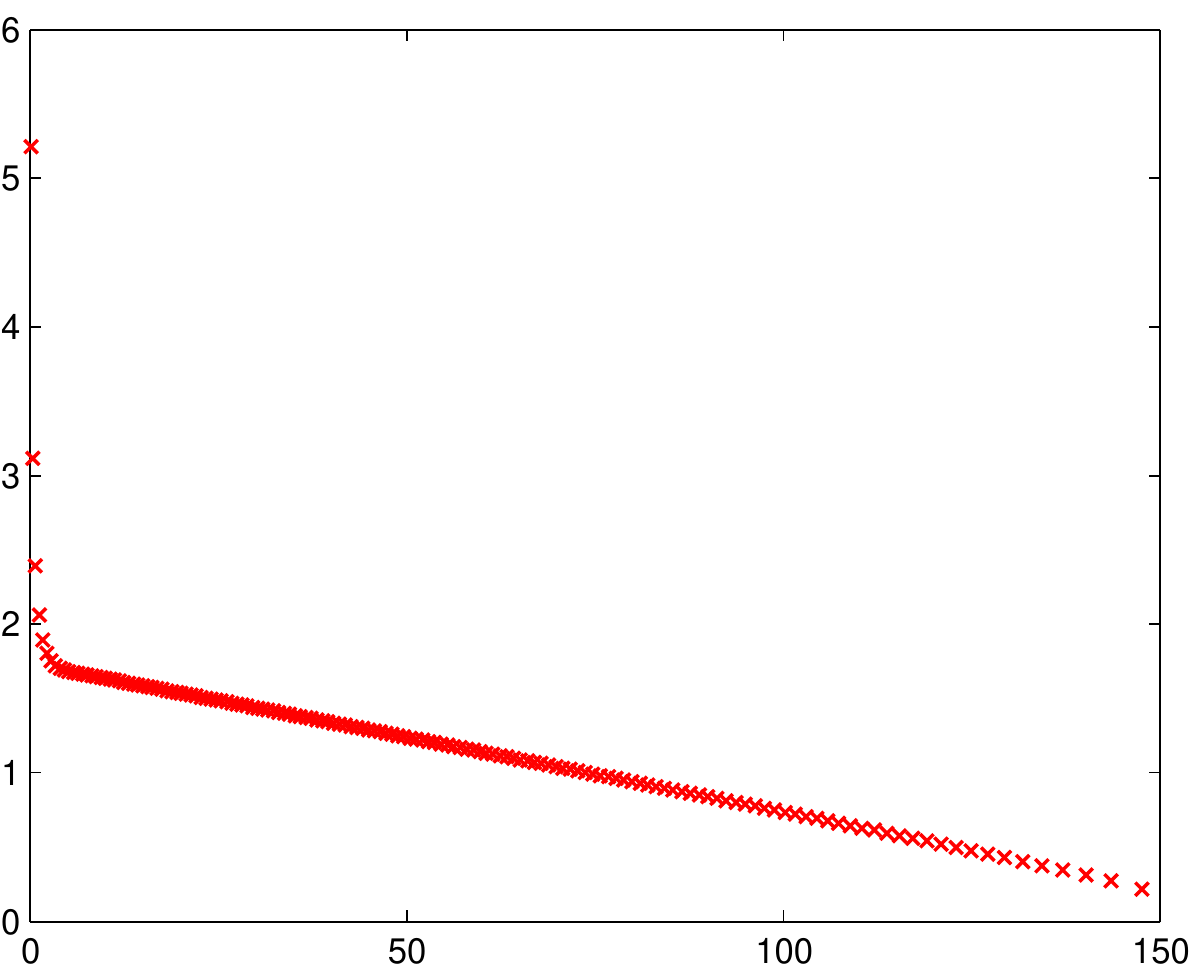} &
\includegraphics[width=1.6in]{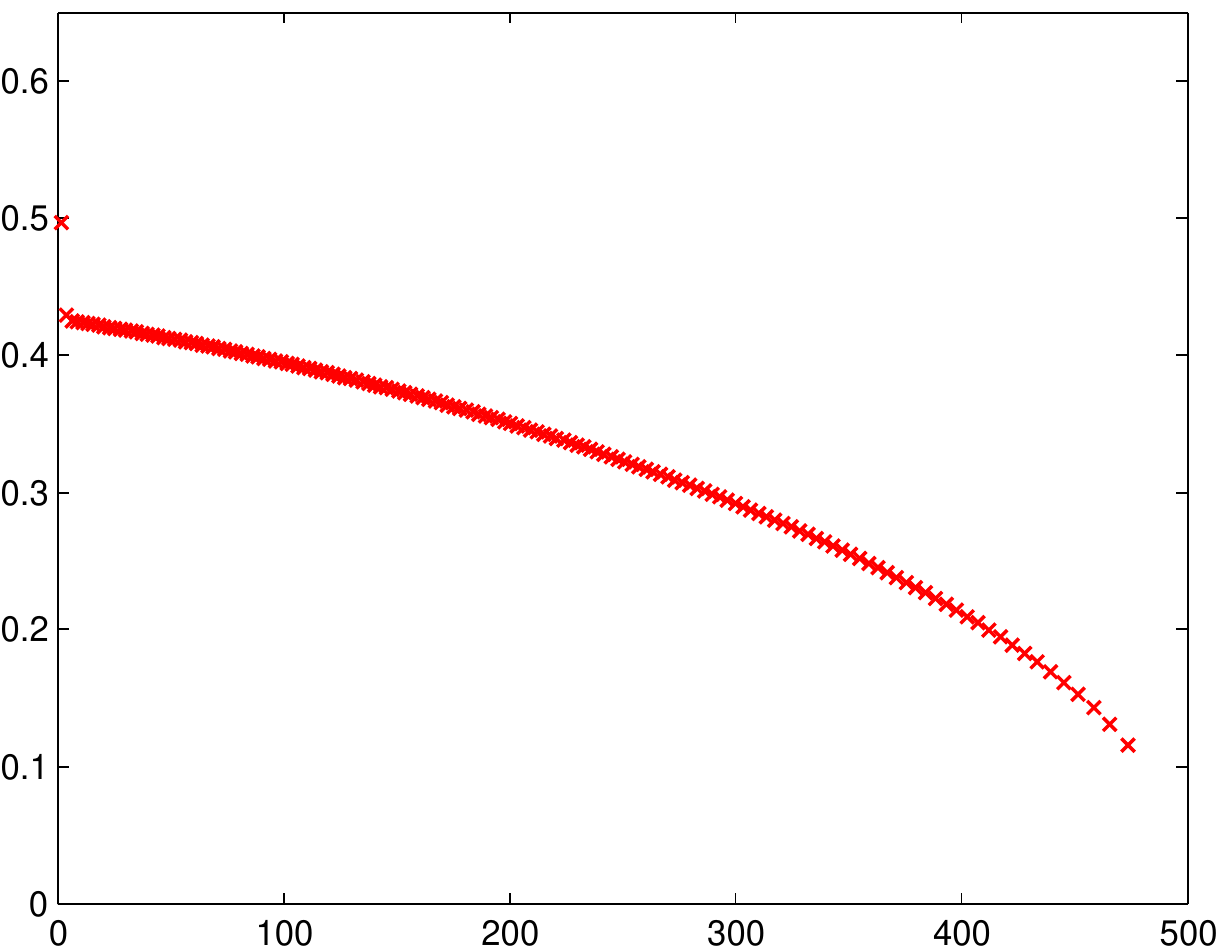}
\end{array}$
\end{center}
\caption{Plots of the discrete dislocation density defined as in (\ref{ddensity})  corresponding to solutions of~(\ref{Eq:1r}) for different values of the applied stress $\sigma$. In all the plots $K=1$, $n=150$, and $h=10$; in (a) $\sigma= 40$; in (b) $\sigma=0.01$; in (c) $\sigma= 0.0005$.}
\label{DDensities}
\end{figure}

Solving the equations (\ref{Eq:1r}) numerically gives a vector $\tilde{x}\in \mathbb{R}^n$; the $i$-th component of the vector represents the  equilibrium position of the $i$-th dislocation wall. Defining the discrete dislocation density as
\begin{equation}\label{ddensity}
\tilde \rho^d(\tilde{x}_i):= \frac{1}{\tilde{x}_i-\tilde{x}_{i-1}}, \qquad i=1,\dots,n,
\end{equation}
we can plot the discrete density as a function of the equilibrium positions, obtaining the plots in Figure \ref{DDensities} for different values of the applied stress $\sigma$.

From Figure \ref{DDensities} we can clearly see that the solution highly depends on the values of the problem parameters. In particular, for large values of the applied stress $\sigma$ the peak of the dislocation density close to zero seems to be dominating, while for smaller values of the stress a rather linear behaviour of the density is prominent at first, and then a decay.
This suggests that depending on the problem parameters, some behaviour of the density seems to be dominant in the bulk (sufficiently far from the boundary of the pile-up region, where the occurrence of \textit{boundary layers} is expected).

This key observation is already present in the recent paper \cite{TRB}, where a careful numerical analysis of the behaviour of the dislocation density in a pile-up is performed, and three different parameter regimes are identified. In particular, the linear dependence of the dislocation density with respect to the positions of the dislocations was observed there.

However, a discrete numerical approach has some drawbacks. First of all, it is computationally expensive when the number $n$ of dislocation walls is large.
But the main drawback is that it is too detailed: knowing the position of each dislocation is definitely not necessary, i.e., an average mesoscopic information in terms of e.g.\ a continuum dislocation density should be preferred.

\section{Discrete energy formulation, scalings, and main results}\label{Main}

Despite its simple formulation, computing the continuum limit of the equilibrium equations~(\ref{Eq:1r}) turns out to be a rather involved problem.
The complexity of the problem is illustrated in Figure~\ref{DDensities}, where the optimal discrete dislocation density (obtained numerically) is plotted for different values of the problem parameters.

\begin{figure}[htbp]
\labellist
\pinlabel {\small Discrete Energy } [t] at 115 376
\pinlabel {\small $\Gamma$-convergence} [t] at 283 390
\pinlabel {\small Continuum Energy } [t] at 466 376
\pinlabel {\small Discrete Equilibrium } [t] at 118 173
\pinlabel {\small Equations } [t] at 125 145
\pinlabel {\small Continuum Equilibrium } [t] at 469 173
\pinlabel {\small Equation } [t] at 465 145
\endlabellist
\begin{center}
\includegraphics[width=3.8in]{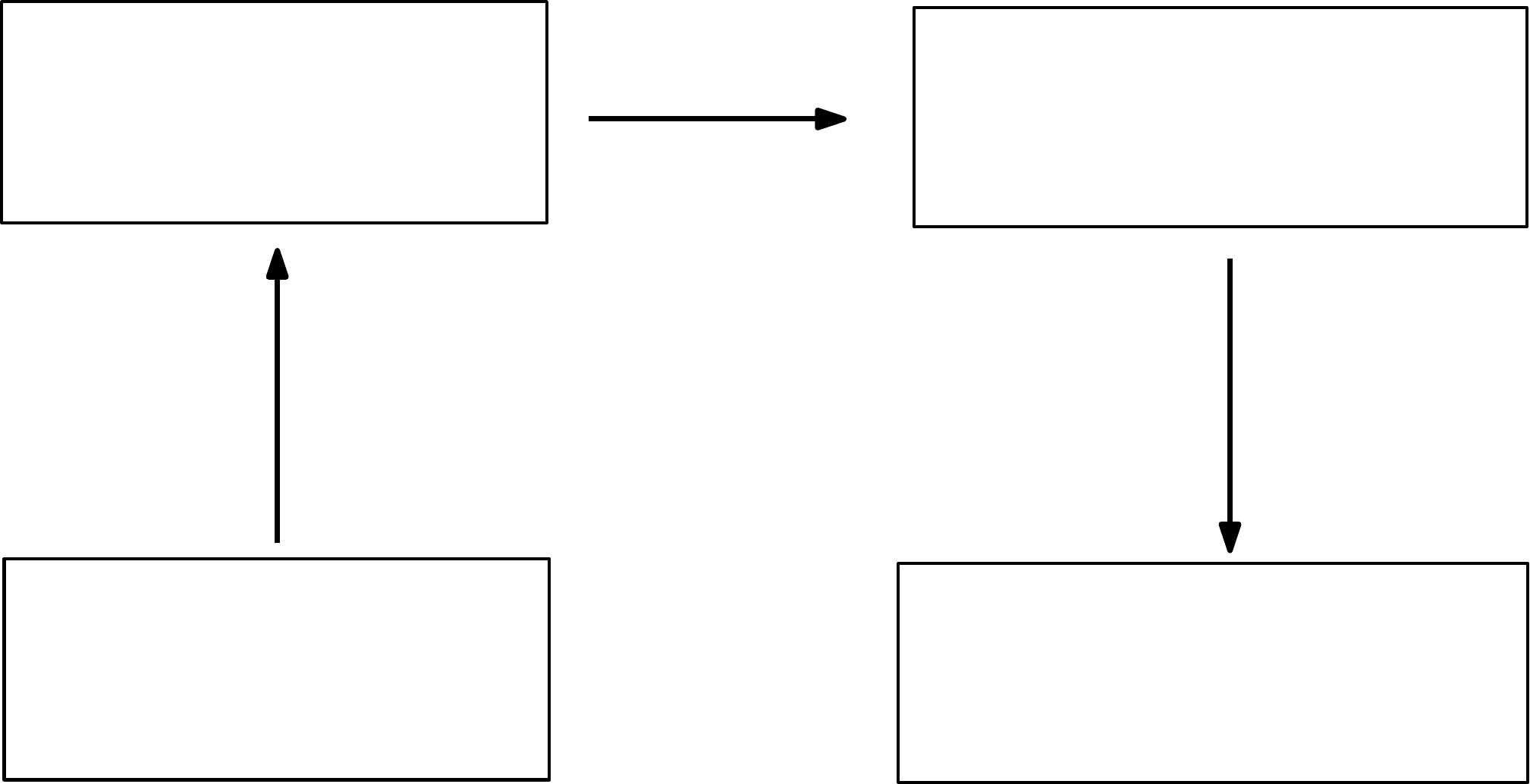}
\end{center}
\caption{Main steps of our upscaling procedure.}
\label{Schema}
\end{figure}
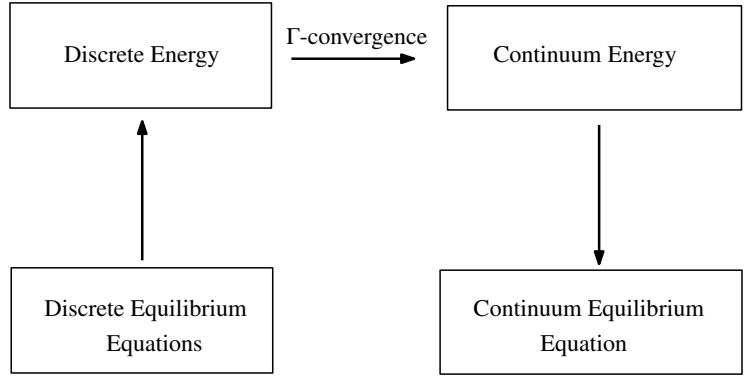

The aim of this paper is to derive a continuum description from the discrete system by rigorous mathematical upscaling.
We do this by exploiting the variational structure of the equilibrium equations (\ref{Eq:1r}); hence we write down the discrete energy whose minimisers satisfy the equilibrium equations (\ref{Eq:1r}). The upscaling procedure will then be carried out on the discrete energy rather than on the
discrete equations, and will produce a continuum energy functional. The method used to perform the discrete-to-continuum upscaling is $\Gamma$-convergence and the details of the rigorous mathematical derivation can be found in \cite{GPPS_Math}. The equilibrium equation associated with the upscaled limit functional is exactly the \textit{upscaled} version of the discrete equations and hence characterises the mesoscopic internal stress as desired. The steps of this procedure are illustrated in Figure \ref{Schema}.

\medskip

We now give an outline of the method and the results. We will discuss these results and put them in perspective, in Sections \ref{Subcritical}-\ref{Comments}.

\subsection{The discrete energy}
\begin{figure}[h]
\labellist
\pinlabel $V(s)$ [tr] at 170 100
\pinlabel $s$ at 345 -10
\pinlabel $\sim -\frac1{\pi^2}\log|s|$ [tl] at 180 220
\pinlabel $\sim\frac2\pi|s|e^{-2\pi|s|}$ [bl] at 220 25
\endlabellist
\begin{center}
\includegraphics[width=2.8in]{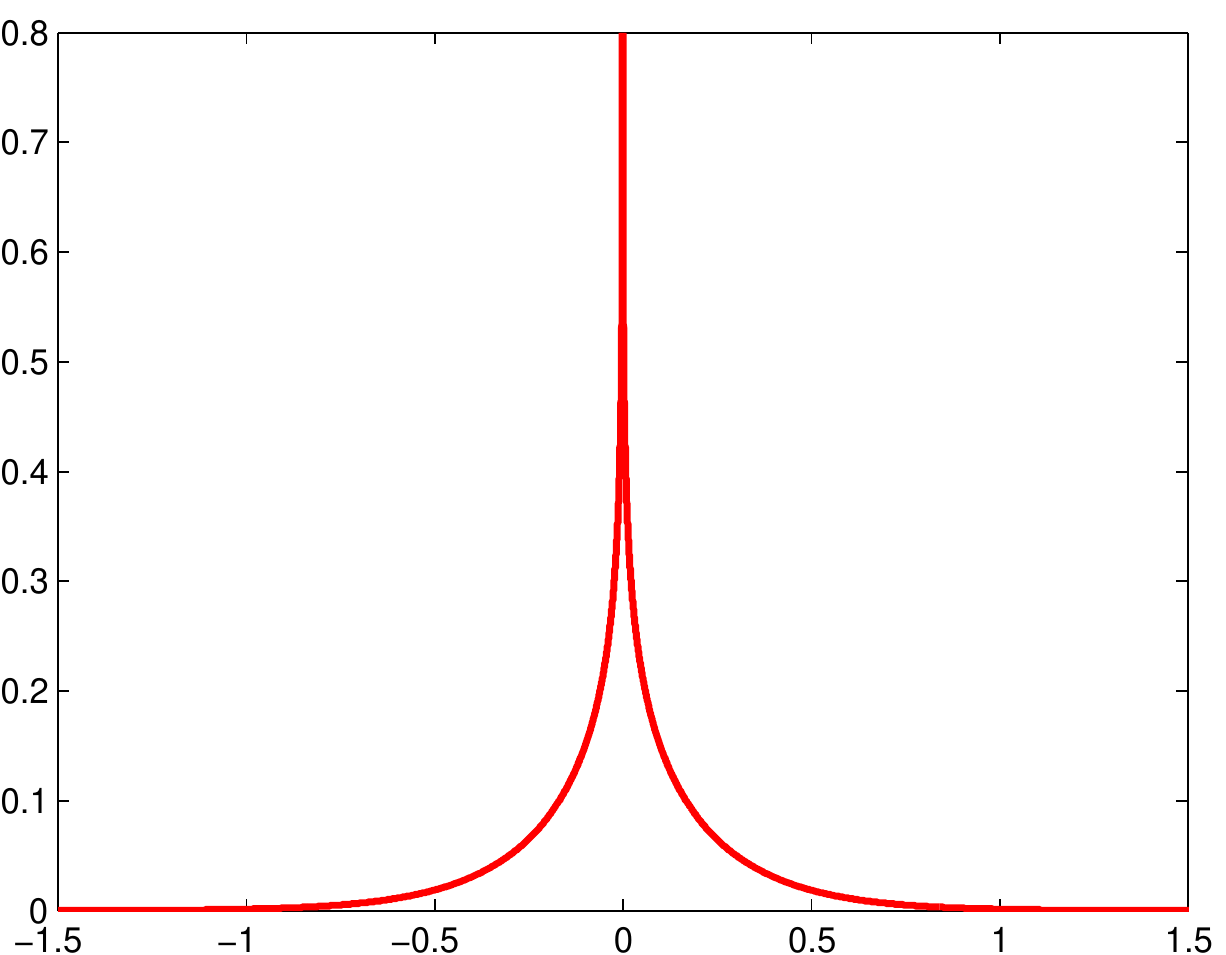}
\end{center}
\caption{The interaction energy  $V$.}
\label{PlotV}
\end{figure}
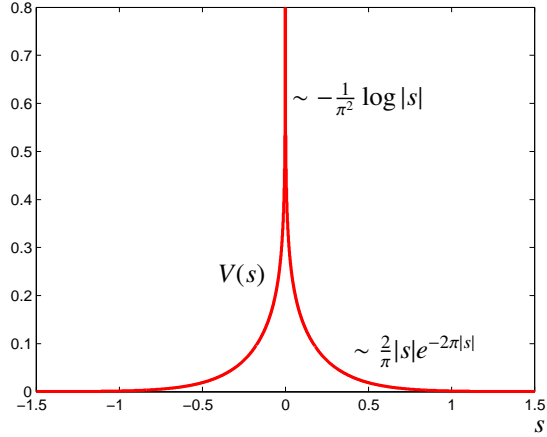

We define the energy associated to the equations (\ref{Eq:1r}). For $\tilde{x}\in \mathbb{R}^n$ let $\tilde{E}_n(\tilde{x})$ be defined as:
\begin{equation}\label{discreteenergy}
\tilde{E}_n(\tilde{x}):=\frac{K}2\sum_{i=0}^n\sum_{\stackrel{j=0}{j\neq i}}^n V\left(\frac{\tilde{x}_i-\tilde{x}_j}{h}\right) +\sigma \sum_{i=0}^n \tilde{x}_i,
\end{equation}
where $V$ is the primitive of $-\varphi$ which decays to zero at infinity, defined as
\begin{equation}\label{defV}
V(s):= \frac{1}{\pi}s\coth \pi s - \frac{1}{\pi^2}\log (2\sinh\pi s),
\end{equation}
and plotted in Figure \ref{PlotV}.
The equations (\ref{Eq:1r}) are equivalent to $\partial_{\tilde{x}_i}\tilde{E}_n=0$.

\medskip

The first term in the energy (\ref{discreteenergy}) penalises configurations where the dislocations are close to one another, since the density $V$ blows up logarithmically in zero, while it favours configurations where dislocations are at large distance from one another, since $V$ decays exponentially at infinity. The second term of the energy, on the contrary, penalises configurations where the dislocations are far from the obstacle at $\tilde{x}=0$.
Using the fact that $V(s)=V(-s)$ we can rewrite the sum in~(\ref{discreteenergy}) as a one-sided sum, as
\begin{equation}\label{DEk}
\tilde{E}_n(\tilde{x})= K\sum_{k=1}^{n}\sum_{j=0}^{n-k} V\left(\frac{\tilde{x}_{j+k}-\tilde{x}_j}h\right) + \sigma\sum_{j=0}^n \tilde{x}_j.
\end{equation}
The form (\ref{DEk}) of the energy will be more convenient for our analysis.

\subsection{Rescaling}

We  convert (\ref{DEk}) into a dimensionless form involving the smallest possible number of independent parameters. This procedure will provide the relevant energy scalings that we  consider in the discrete-to-continuum derivation.

We first non-dimensionalize the dislocation positions $\tilde x_i$ by a length scale $\ell_n$, by introducing the new dimensionless variables $x_i$,
\[
x_i := \frac {\tilde x_i}{\ell_n}.
\]
At this stage the length scale $\ell_n$ is undetermined; we will use our freedom of choice of $\ell_n$ below to capture the behaviour of the pile-up. This is akin to choosing the appropriate degree of magnification in an experimental setup. The choices that we make below will make $\ell_n$  of the same order as the pile-up length.

Given a choice of $\ell_n$, the \textit{aspect ratio} $\alpha_n$ is the ratio between the average dislocation distance $\Delta \tilde{x}\sim \frac{\ell_n}{n}$ and the vertical spacing $h$ between consecutive dislocations in the same wall:
\[
\alpha_n:=\frac{\ell_n}{nh}.
\]
We now rewrite the energy (\ref{DEk}) in a dimensionless form, in terms of the new variable $x$:
\begin{equation}\label{Energy2}
E_n(x):= \frac{\tilde{E}_n(n h\alpha_n x)}{n^2\sigma h\alpha_n}= \frac{K}{n^2\sigma h \alpha_n}\sum_{k=1}^{n}\sum_{j=0}^{n-k} V\left(n\alpha_n(x_{j+k}-x_j)\right) + \frac{1}{n}\sum_{j=0}^n x_j.
\end{equation}
Since $E_n$ will be of order one, the dimensional energy will scale as $n^2\sigma h \alpha_n$.

The aspect ratio $\alpha_n$ is a local description of the average arrangement of the dislocations within the pile-up region, whereas $\ell_n$ gives a more global information on the dislocation distribution. Therefore there are two ratios of length scales: $\ell_n/h$, which compares the total length of the pile-up with the vertical spacing $h$, and $\alpha_n$, which compares the spacing between the walls and the spacing within the walls.

\begin{figure}[h]
\labellist
\pinlabel {\small Inaccessible} [t] at 120 450
\pinlabel {\small (1)} [t] at 150 360
\pinlabel \small $\alpha_n$ [t] at 45 500
\pinlabel  $\frac{\ell_n}{h}=n\alpha_n$ [t] at 270 305
\pinlabel \small $1$ [t] at 48 415
\pinlabel \small $1$ [t] at 170 305
\pinlabel {\small (2)} [t] at 30 395
\pinlabel {\small (5)} [t] at 205 430
\pinlabel {\small (3)} [t] at 205 365
\pinlabel {\small (4)} [t] at 170 505
\pinlabel {\small $\ell_n$} [t] at 400 100
\pinlabel {\small $h$} [t] at 480 150
\pinlabel {\small $\ell_n$} [t] at 400 350
\pinlabel {\small $h$} [t] at 475 425
\pinlabel {\small $\ell_n$} [t] at 130 90
\pinlabel {\small $h$} [t] at 200 190
\endlabellist
\begin{center}
\includegraphics[width=4in]{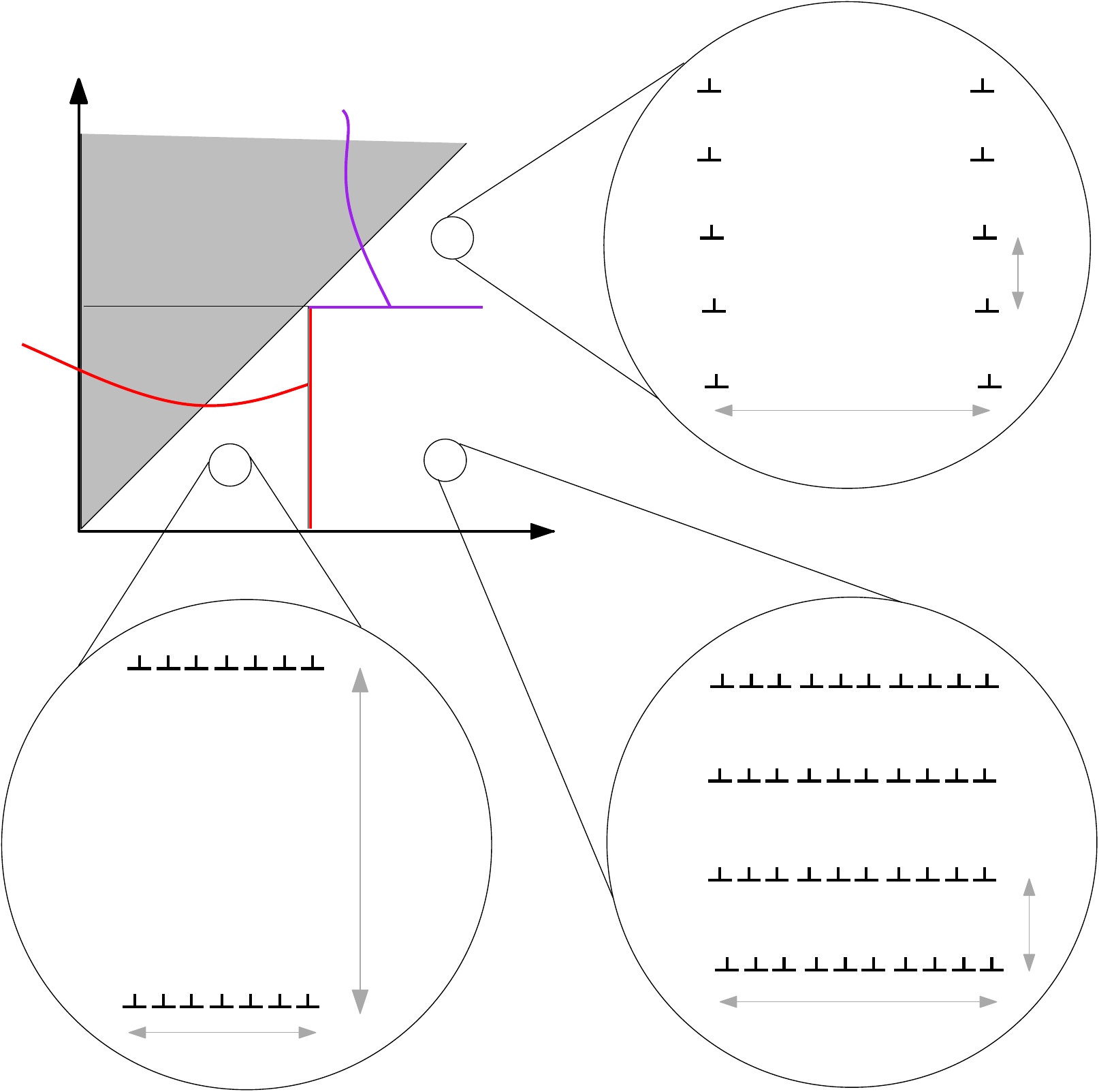}
\end{center}
\caption{Scalings (1)-(5) in the parameter space $\left(\frac{\ell_n}{h},\alpha_n\right)$ and some of the corresponding configurations.}
\label{PUregions}
\end{figure}
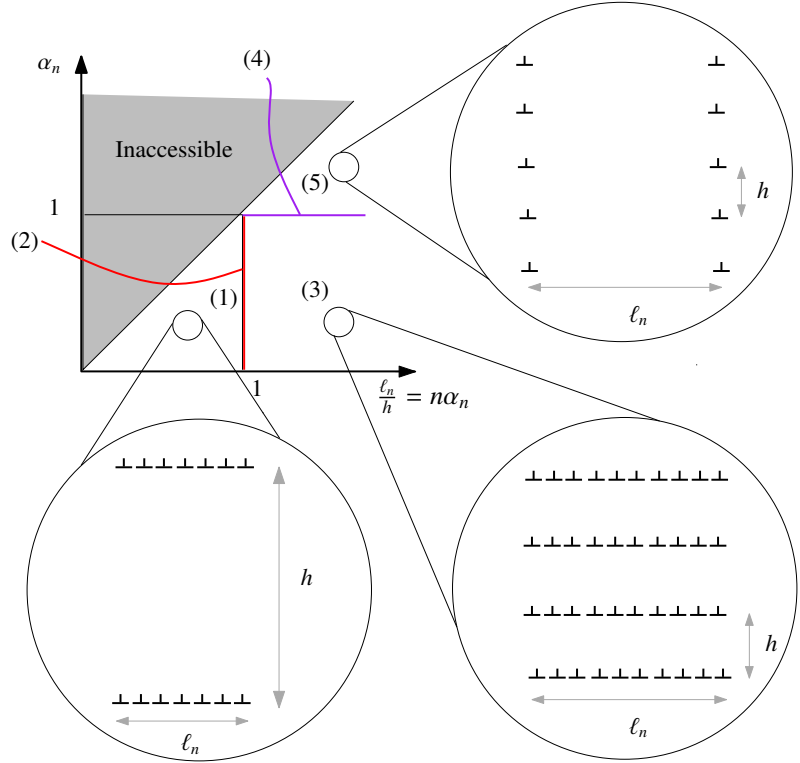

The possible relations between these two (clearly non-independent) parameters is illustrated in Figure \ref{PUregions}.
The region marked as inaccessible corresponds to the regime $\alpha_n > \frac{\ell_n}{h}$, which is not compatible with the definition of $\alpha_n$ since it would require that the length of the pile-up region is smaller than the spacing between the walls.

The scalings of $\alpha_n$ illustrated in Figure \ref{PUregions} and that we will describe are listed below:

\bigskip

$
\label{scalings-alpha}
\left.
\begin{array}{ll}
\medskip
(1)\,\, \mbox{Subcritical regime: } \alpha_n\ll \frac{1}{n};\\
\medskip
(2)\,\, \mbox{First critical regime: } \alpha_n\sim \frac{1}{n};\\
\end{array}
\,\,\qquad\right\}
\quad \frac{\ell_n}{h} \lesssim 1
$

$
\left.
\begin{array}{lll}
\medskip
(3)\,\, \mbox{Intermediate regime: } \frac1n \ll \alpha_n\ll 1;\\
\medskip
(4)\,\, \mbox{Second critical regime: } \alpha_n\sim 1;\\
\medskip
(5)\,\, \mbox{Supercritical regime: } \alpha_n\ll 1.
\end{array}
\right\}
\quad \frac{\ell_n}{h}\gg 1
$

\subsection{Determining the scaling by a trial configuration}

In the companion paper~\cite{GPPS_Math} we prove that, given certain choices of $\ell_n$ for each of the cases $(1)$-$(5)$ above, the rescaled energy $E_n$ converges as $n\to\infty$. Here we do not give the details of this proof; instead we only specify the choices of $\ell_n$, but motivate, heuristically, why these are the right ones, and how the resulting limit expressions for the energy and the internal stress arise.

A convenient method to determine the scaling of $\ell_n$ is by requiring that both terms of the energy (\ref{Energy2}) are of the same order. As it turns out, the correct scaling can be determined by only considering the special, uniformly spaced configuration $x_i := i/n$. In this case, the second term in~(\ref{Energy2}) is always of order one, since
\[
\frac1n \sum_{j=0}^n x_j \to \frac12 \quad\text{as } n\to\infty,
\]
and therefore we will require that the choice of $\ell_n$ makes the first term in~(\ref{Energy2}) also of order one.

By this analysis we are able to relate $\ell_n$ and
$\alpha_n$ to the parameters $h$, $\sigma$, $K$ and the number~$n$ of the dislocation walls. In particular, we  relate $\ell_n$ and
$\alpha_n$ to the dimensionless parameter~$\beta_n$ defined as
\begin{equation}\label{betan}
\beta_n:= \sqrt{\frac{K}{n\sigma h}},
\end{equation}
introduced in \cite{TRB}. In mechanical terms, $\beta_n$ is strongly related to the forcing term $\sigma$ and is small for large forcing and large for small forcing.

Note that this implies that $K$, $\sigma$, and $h$ may depend on $n$, and in many cases this  makes sense. For instance, if the number of dislocations $n$ increases because the density of Frank-Read sources increases, then one would expect that the average slip plane spacing $h$ should simultaneously decrease with $n$.

Performing this matching on the trial configuration of uniformly spaced dislocations, as we do in Sections~\ref{Subcritical}--\ref{Supercritical} below, we find that both terms in~(\ref{Energy2}) have the same scaling provided we make the following choices:
\begin{enumerate}
\item[(a)]
If $\beta_n\ll 1/n$, then choose
\[
\ell_n = \frac{Kn}\sigma,\qquad  \alpha_n = \frac K{\sigma h}.
\]
Then $\alpha_n\ll 1/n$, and therefore this is the Subcritical regime.
\item[(b)]
If $\beta_n \sim 1/n$, $1/n \ll \beta_n \ll1$, or $\beta_n\sim 1$, then choose
\[
\ell_n = \sqrt{\frac{Knh}\sigma}, \qquad \alpha_n = \beta_n =\sqrt{\frac K{nh\sigma}}.
\]
Since $\alpha_n= \beta_n$, these three cases correspond to the First critical regime, the Intermediate regime, and the Second critical regime, respectively.
\item[(c)]
If $\beta_n\gg 1$, then choose
\[
\ell_n = \frac1{2\pi} n^2 h\sigma \ln \Big(\frac{2K}{nh\sigma}\Big),\qquad
\alpha_n = \frac1{2\pi} n\sigma \ln \Big(\frac{2K}{nh\sigma}\Big).
\]
Then $\alpha_n\gg 1$, and this is the Supercritical regime.
\end{enumerate}
In this way we are able to determine in which regime the system is  by considering only $\beta_n$, which is an explicit function of the constants, and independent of the solution.

\subsection{Internal stress}

The most convenient concept to connect the discrete (finite $n$) situation with its infinite-$n$ limit is the \emph{density} $\tilde \rho^d$ introduced in (\ref{ddensity}). As $n$ becomes large, this discrete density approximates a continuous density $\tilde \rho$.
For four of the five regimes we derive expressions for the internal stress in terms of this limiting, continuous  density $\tilde \rho$:

\begin{align}
&(1)\,\, \mbox{Subcritical regime: }& &\sigma^{(1)}_{\textrm {\scriptsize{int}}}(\tilde x)  =  \frac{K}{\pi^2}\int_0^\infty \log{|\tilde x-\tilde y|}\,\partial_{\tilde y}\tilde \rho(\tilde y)\, d\tilde y;
\label{intstress1dimensional}\\
&(2)\,\, \mbox{First critical regime (assuming $\beta_n=1/n$): }& &\sigma^{(2)}_{\textrm{int}}(\tilde x) = - K\int_0^\infty V\Bigl(\frac{\tilde x-\tilde y}h\Bigr) \, \partial_{\tilde y}\tilde \rho(\tilde y)\, d\tilde y ;\label{intstress2dimensional}\\
\bigskip
&(3)\,\, \mbox{Intermediate regime: } & &\sigma^{(3)}_{\textrm{int}}(\tilde x) = - \frac{Kh}{3\pi} \partial_{\tilde x}  \tilde \rho;\label{intstress3dimensional}\\
\bigskip
&(4)\,\, \mbox{Second critical regime (assuming $\beta_n=1$):}& &\sigma^{(4)}_{\textrm{int}}(\tilde x) = \frac {K}{h^2} \frac{1}{\tilde \rho^3}\, \varphi'_{\textrm{eff}}\left(\frac{1}{h\tilde \rho}\right)\, \partial_{\tilde x}\tilde \rho.\label{intstress4dimensional}
\end{align}
Here
\[
\varphi_{\mathrm{eff}}(t) := \sum_{k=1}^\infty k\varphi(kt) = - \sum_{k=1}^\infty kV'(kt),
\]
and $\varphi'_{\mathrm{eff}}$ has the asymptotic behaviour (see Section~\ref{SecondCritical})
\[
\varphi'_{\mathrm{eff}}(t) \sim \begin{cases}
\displaystyle -\frac1{3\pi t^3} & \text{as }t\to0\\[3\jot]
\displaystyle -8\pi te^{-2\pi t} & \text{as }t\to\infty.
\end{cases}
\]
Note that the argument of $\varphi'_{\mathrm{eff}}$, $1/h\tilde\rho$, is an aspect ratio, since it is the ratio of the spacing $(\tilde x_i-\tilde x_{i-1})$ to~$h$.

The fifth, Supercritical, regime, has no meaningful internal stress, since it is too degenerate (see Section~\ref{Supercritical}).

\subsection{Pile-up behaviour}

For each of the four cases above, the equilibrium density $\tilde \rho$ of a set of dislocation walls that is pushed against the obstacle at $\tilde x = 0$ solves the equation
\begin{equation}
\label{eq:sigmas}
\sigma_{\mathrm{int}}(\tilde x) - \sigma = 0.
\end{equation}
\begin{enumerate}
\item For the Subcritical regime,  equation~(\ref{eq:sigmas}) has been solved by Head and Louat~\cite{HL}:
\[
\tilde \rho(\tilde x) = \frac\sigma K \sqrt{\frac{\ell-\tilde x}{\tilde x}}, \qquad
\text{with } \ell = \frac{2nK}{\pi^2\sigma}.
\]
\item For the First critical regime, equation~(\ref{eq:sigmas}) has no explicit solution, to our knowledge.
\item In the Intermediate regime, the equilibrium density is given by
\[
\tilde \rho(\tilde x) = \frac{2n}{\ell^2} (\ell-\tilde x), \qquad
\text{with } \ell = \sqrt{\frac{2nKh}{3\pi\sigma}}.
\]
\item In the Second critical regime, again we believe that equation~(\ref{eq:sigmas}) has no known explicit solution.
\item Finally, in the Supercritical regime, the equilibrium density $\tilde \rho$ is constant:
\[
\tilde\rho (\tilde x) = \begin{cases}
\ell^{-1}
& \text{if }\tilde x \leq \ell, \\
0 & \text{otherwise,}
\end{cases}
\qquad \text{with } \ell = \frac1{2\pi} n^2 h\sigma \ln \Big(\frac{2K}{nh\sigma}\Big).
\]
\end{enumerate}
In the three explicit cases above, the parameter $\ell$ is exactly the pile-up length. Note that $\ell$ scales the same as the  scaling parameter $\ell_n$, but the two differ by a numerical constant. In Figures~\ref{DD1}--\ref{Fig:constant} the continuum pile-up densities are compared to their discrete counterparts.

\bigskip

We  now consider each of the five regimes separately in Sections \ref{Subcritical}--\ref{Supercritical}; in every case we validate the continuum upscaled model by comparing its predictions with the solution of the discrete model. A comparison between the internal stresses $\sigma_{\textrm{int}}^{(k)}$ above (for $k=1-5$) with other well-known models in the engineering literature will be the subject of Section \ref{Comparison}. Section \ref{Comments} contains a number of interesting comments on the the results of this paper.

\section{Subcritical regime $\alpha_n\ll\frac1n$: derivation of the upscaled internal stress}\label{Subcritical}

In this section we consider the extreme case $\alpha_n\ll\frac1n$, or equivalently, $\frac{\ell_n}{h}\ll 1$, i.e., the case where the length $\ell_n$ of the pile-up region is much smaller
than the in-wall dislocation spacing. In this case the in-plane interaction is much stronger than the in-wall interaction, i.e., the equivalent continuum
formulation will not \textit{sense} the walls and the result will correspond to the case of a single slip plane rather than infinite walls of dislocations.

This scaling regime has been analysed numerically (in its discrete formulation) in the paper~\cite{TRB}, and corresponds to the so-called \textit{close region} in their terminology.

\subsection{Heuristics for the scaling of the discrete energy}
As described above, we will identify  $\ell_n$ and $\alpha_n$  and the corresponding rescaling of the energy in~(\ref{Energy2}) by requiring that the two terms of the discrete energy $E_n$ are of the same order and bounded when calculated for the uniformly-spaced configuration.

Since the second term in~(\ref{Energy2}) is of order one, the only condition we are left to impose  is that the first term of the energy is of order one, i.e.,
\begin{equation}\label{C:2}
\frac{K}{n\sigma h\alpha_n}\sum_{k=1}^{n}V\left(\alpha_nk\right) \sim 1.
\end{equation}
Since for this rescaling $n\alpha_n \ll 1$, the argument of $V$ in (\ref{C:2}) is small, as $\alpha_n k\leq \alpha_n n \ll1$. Therefore we substitute for $V$ in (\ref{C:2}) its asymptotic expansion close to zero,
\[
V(s)\sim \frac{1-\log 2\pi s}{\pi^2}, \quad \mbox{for } s>0.
\]
Using this approximation in (\ref{C:2}) we find that we should require that
\begin{equation}\label{C:2a}
\frac{K}{n\sigma h\alpha_n}\sum_{k=1}^{n}\left[\frac{1- \log(2\pi \alpha_n k)}{\pi^2}\right] \sim 1.
\end{equation}
As it stands, this condition would give rise to a wrong choice of $\ell_n$ and $\alpha_n$. This is because  the expression between brackets is dominated by a constant term, and this constant term is irrelevant for the equilibrium equations (since it vanishes upon differentiation). We therefore introduce the renormalized energy density $\hat{V}_n(t):= V(t) +\frac{-1 + \log(2\pi n \alpha_n)}{\pi^2}$ in which we have subtracted  this large constant.  In terms of $\hat{V}_n$ the bound (\ref{C:2a}) on the energy becomes
$$
- \frac{1}{\pi^2}\frac{K}{n\sigma h\alpha_n}\sum_{k=1}^{n}\log\left(\frac{k}{n}\right) \sim 1.
$$
Simple computations show that
\begin{align}\label{RSlog}
\frac{K}{n\sigma h\alpha_n}\sum_{k=1}^{n}\log\left(\frac{k}{n}\right) &=
\frac{K}{\sigma h\alpha_n}\left[\frac1n\sum_{k=1}^{n}\log\left(\frac{k}{n}\right)\right] \sim \frac{K}{\sigma h\alpha_n},
\end{align}
since the term in square brackets is the Riemann sum for the integral $\int_0^1\log(t)\,dt = 1$. Hence the bound on the energy reduces to the requirement that $\frac{K}{\sigma h\alpha_n} \sim 1$. This bound provides the following expressions for the aspect ratio $\alpha_n$ and for the length $\ell_n$ of the pile-up region in terms of the parameters $\sigma$, $h$, $K$ and $n$:
\begin{equation}\label{al1}
\alpha_n^{(1)} \sim \frac{K}{\sigma h}; \quad \ell^{(1)}_n \sim \frac{Kn}{\sigma},
\end{equation}
or, in terms of the dimensionless parameter $\beta_n$ defined in (\ref{betan}),
\begin{equation}\label{ab1}
\alpha^{(1)}_n \sim  n\beta_n^2; \quad \ell^{(1)}_n \sim n^2\beta_n^2 h.
\end{equation}
We note that, by (\ref{ab1}), the scaling regime $\alpha_n\ll\frac{1}{n}$ can be equivalently formulated in terms of~$\beta_n$ and corresponds to
$\beta_n\ll\frac{1}{n}$.
From \eqref{Energy2}, using \eqref{al1} and \eqref{ab1}, we obtain the scaling of the discrete energy in the subcritical regime, which reads

\begin{align}\label{E1:d}
E_n^{(1)}(x) :=
\frac{1}{n^2}\sum_{k=1}^{n}\sum_{j=0}^{n-k} \hat{V}_n\left(n^2\beta_n^2(x_{j+k}-x_j)\right) + \frac1n\sum_{j=0}^n x_j.
\end{align}
Note that $\hat{V}_n$ can be expressed in terms of $\beta_n$ as $\hat{V}_n(t)= V(t) +\frac{-1 + \log(2\pi n^2\beta_n^2)}{\pi^2}$.

\subsection{Continuum limit: derivation of the internal stress}

In \cite{GPPS_Math} we rigorously derived a continuum energy in terms of a continuum dislocation density $\rho$, starting from the discrete energy \eqref{E1:d}, in the limit $n\to \infty$.
Here we give a hint of the main idea of the proof, and we refer to \cite{GPPS_Math} for the detailed derivation.

Starting from the discrete energy \eqref{E1:d} we define $\rho_n:= \frac{1}{n}\sum_{i=1}^n\delta_{x_i}$, where $\delta_{x_i}$ is the Dirac delta function  localised at $x_i$ (and is zero everywhere except at $x_i$). The measure $\rho_n$ describes the distribution of the walls and is approximately a dimensionless version of the discrete density $\tilde \rho^d$ introduced in \eqref{ddensity}. In terms of $\rho_n$ we can rewrite the sums in \eqref{E1:d} as integrals, namely
\begin{equation}\label{E1:dna}
E_n^{(1)}(x)= \frac{1}{2}\int_0^\infty\int_0^\infty \hat{V}_n\left(n^2\beta_n^2(x-y)\right)\rho_n(x)\rho_n(y)dx dy + \int_0^\infty x \rho_n(x)dx.
\end{equation}
In the previous formula we wrote the integral on $(0,\infty)$ since the wall positions $x_i$ are the unknowns we are solving for, which are in $(0,\infty)$.
We note that in terms of $\rho_n$ the discrete energy looks like a continuum functional. By the definition of $\hat{V}_n$ and since $n\beta_n\ll 1$ in this regime, we have that $\hat{V}_n\left(n^2\beta_n^2s\right)\sim -\frac{1}{\pi^2}\log|s|$. Substituting this expression in \eqref{E1:dna} we have
\begin{equation}\label{E1:dn}
E_n^{(1)}(x)\simeq -\frac{1}{2\pi^2}\int_0^\infty\int_0^\infty \log|x-y|\,\rho_n(x)\rho_n(y)dx dy + \int_0^\infty x \rho_n(x)dx.
\end{equation}
It turns out that for a large number of walls, i.e., as $n\to \infty$, $\rho_n$ converges to a continuum density~$\rho$ and the discrete energy $E^{(1)}_n$ converges to the continuum functional $E^{(1)}$ defined as:
\begin{equation}\label{E1:c}
E^{(1)}(\rho):= -\frac{1}{2\pi^2}\int_0^{\infty}\int_0^{\infty} \log|x-y|\,\rho(x)\rho(y)dx dy  + \int_0^\infty x\rho(x)dx.
\end{equation}
This is the top arrow in Figure~\ref{Schema}.

According to the scheme in Figure \ref{Schema}, it now remains to compute the Euler-Lagrange equation associated with the continuum functional $E^{(1)}$. In this case, the continuum equilibrium equation is an integral equation of the following form:

\begin{equation}\label{EL:1a}
-\frac{1}{\pi^2}\int_0^{\infty}\frac{\rho(y)}{|x-y|}dy + 1 = 0,
\end{equation}
for every $x\in (0,\infty)$. This equation is to be interpreted as a Cauchy principal-value integral (see e.g.~\cite[Ch.~2]{Musk}). Alternatively, we can rewrite the equation \eqref{EL:1a} in the compact form
\begin{equation}\label{EL:1}
\frac{1}{\pi^2}\log\ast\,\partial_x \rho + 1 =0.
\end{equation}
We recall that our starting point was a discrete system of the form $\sigma_i^{\textrm{int}} - \sigma = 0$. Therefore the continuum equation \eqref{EL:1} can be interpreted as the dimensionless form of $\sigma^{(1)}_{\textrm{int}} - \sigma = 0$. Hence, the expression for the dimensionless internal stress we obtained is
\begin{equation}\label{intstress1}
\sigma^{(1)}_{\textrm{int}} = -\frac{1}{\pi^2}\log\ast\,\partial_x \rho.
\end{equation}
After returning to dimensional variables we find~\eqref{intstress1dimensional}.

\subsection{The Eshelby-Frank-Nabarro model (EFN): the single slip plane case}
Since in this regime the length of the pile-up region $\ell_n$ is much smaller than the vertical spacing $h$, the in-wall interaction
is much weaker than the in-plane interaction. Therefore it is natural to compare our model with the EFN model, formulated in the case of a single slip plane.

In the EFN model the equilibrium positions of $n$ dislocations (as opposed to walls) in a single slip plane (under an applied stress $\sigma$ that pushes them towards the barrier) are described by equations similar to \eqref{Eq:1}. In that case, the interaction potential is $\psi(t):=1/\pi^2 t$.

\begin{figure}[htbp]
\labellist
\pinlabel {\small Dislocation wall positions} at 200 -8
\pinlabel \small $\rho^d(x_i)$ at -30 130
\pinlabel {\scriptsize EFN} at 325 257
\pinlabel {\scriptsize (1)} at 325 244
\endlabellist
\begin{center}
\includegraphics[width=3.6in]{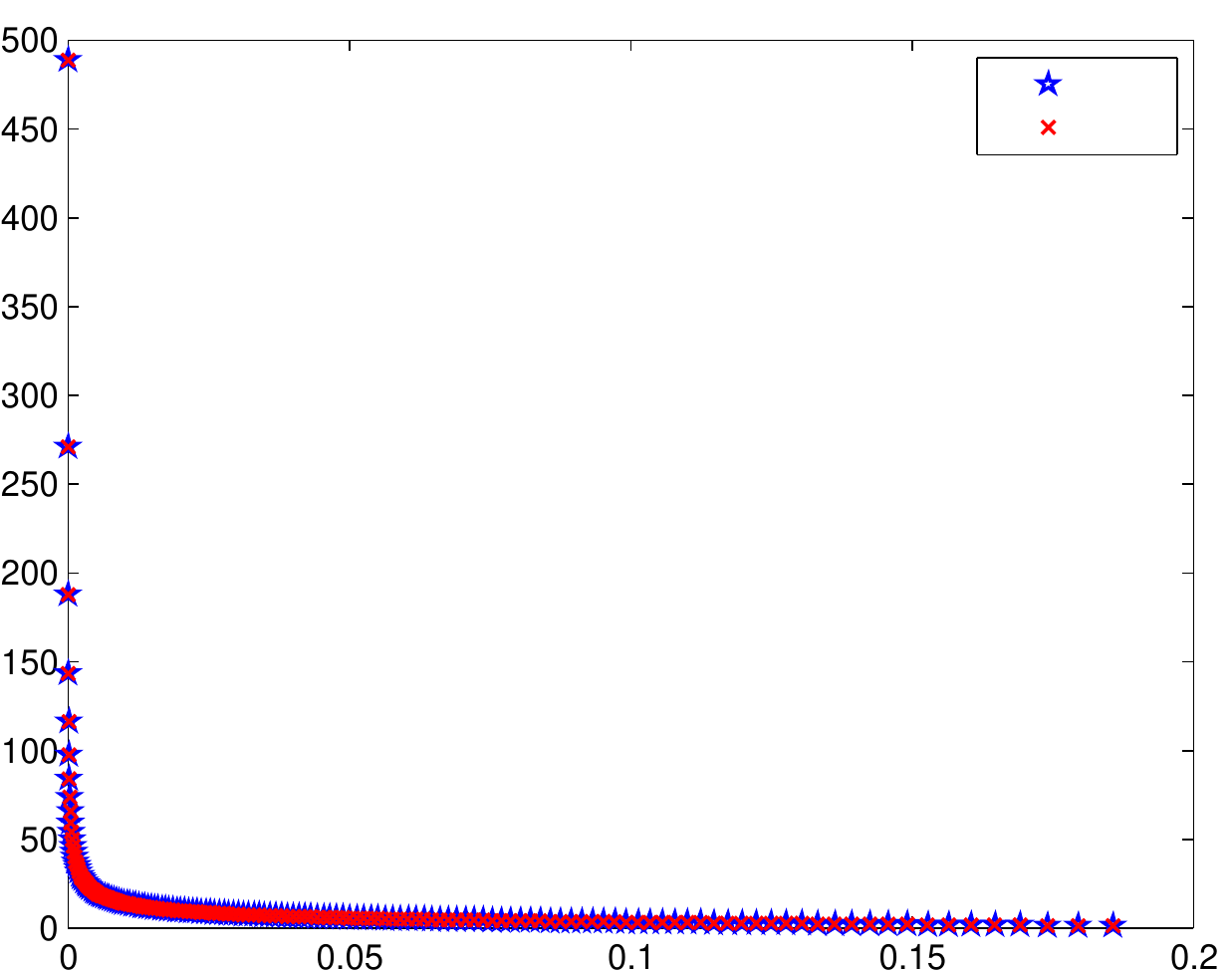}
\end{center}
\caption{Comparison between the optimal discrete density in the EFN model and in~\eqref{E1:d}, in the subcritical scaling regime (1). Here $n=150$ and $\beta_n = \frac{4}{n\sqrt{n}} \left(\ll\frac{1}{n}\right)$. The stars and crosses completely coincide.}
\label{EFNvsOurs}
\end{figure}
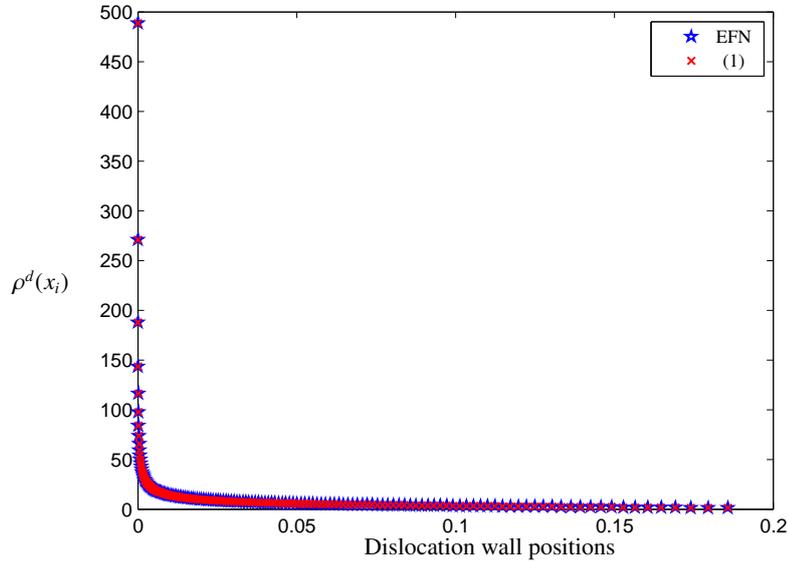

In Figure \ref{EFNvsOurs} we show the comparison between the EFN model and our model. We plot the discrete density $\rho^d$ corresponding to the minimisers $x\in \mathbb{R}^n$ of the energy \eqref{E1:d}, for $\beta_n \ll \frac{1}{n}$, and of the EFN energy, where $\rho^d(x_i):=1/(x_i - x_{i-1})$. The discrete densities corresponding to the two models show perfect agreement.

\subsection{Comparison Discrete vs Continuum}
\begin{figure}[htbp]
\labellist
\pinlabel {\small Dislocation walls positions} [b] at 180 -18
\pinlabel \small $\rho^d$ at 328 246.5
\pinlabel \small $\rho$ at 326 256
\endlabellist
\begin{center}
\includegraphics[width=3.6in]{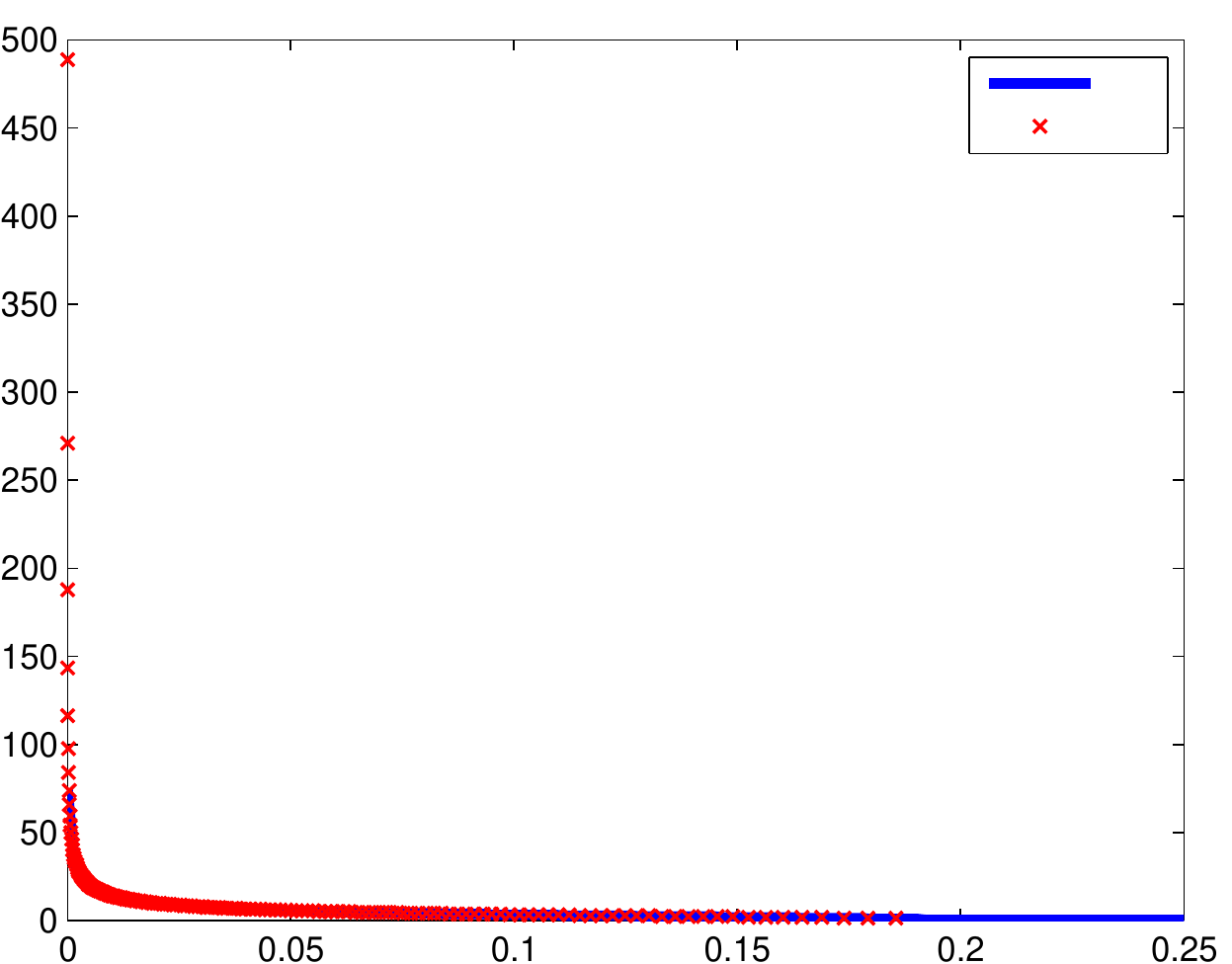}
\end{center}
\caption{Comparison of the discrete and continuous pile-ups, for $n=150$ and $\beta_n = 6/(n\sqrt n)$. The continuum dimensionless density $\rho$ minimizes $E^{(1)}$; the discrete dimensionless density $\rho^d$ is calculated by minimizing $E_n^{(1)}$ over all $x_i$, and defining $\rho^d(x_i):= 1/(x_i-x_{i-1})$. Both are normalized such that the total density equals $1$.}
\label{DD1}
\end{figure}

In this section we numerically compare  the discrete density obtained by minimising the discrete energy \eqref{E1:d} for $\beta_n\ll \frac1n$ and for large $n$ with the solution of the continuum equation \eqref{EL:1}. The agreement is shown in Figure \ref{DD1}.

We note that the continuum equation \eqref{EL:1} can be solved by means of the Hilbert transform, and it has a closed-form solution, found in~\cite{HL} (see also~\cite{Vos}).

\section{First critical regime $\alpha_n\sim\frac1n$: derivation of the upscaled internal stress}\label{FirstCritical}
This regime corresponds to $\ell_n \sim h$, i.e., to the case where the length of the pile-up region is of the same order as the vertical dislocation spacing $h$ and consequently, for large $n$, the horizontal wall-wall spacing is much smaller than $h$. Unlike the previous case, we cannot expect that the vertical interactions can be neglected and that the stress exerted by a wall is equivalent to the stress generated by a single dislocation---at least, not quantitatively. Qualitatively, though, the optimal dislocation density exhibits a sharp increase close to the obstacle, as in the previous case, and a fast decay at infinity.

\subsection{Heuristics for the scaling of the discrete energy}
To find the scaling of the discrete energy $E_n$ in \eqref{Energy2} that guarantees that both terms are bounded and of the same order (and the corresponding expressions for $\alpha_n$ and $\ell_n$) we use as a test configuration, as in the previous scaling, the uniform distribution $x_i=\frac{i}{n}$.
Since with this choice the linear term in the energy is of order one, the bound on the energy reduces to a bound on the first term only, i.e., proceeding as in \eqref{RSlog},
\begin{align}\label{Riemann2}
1\sim \frac{K}{n^2\sigma h \alpha_n}\sum_{k=1}^{n}\sum_{j=0}^{n-k} V\left(\alpha_nk\right) \sim \frac{K}{n\sigma h \alpha_n^2} \left[\alpha_n\sum_{k=1}^{n} V\left(\alpha_nk\right)\right].
\end{align}
Since for this scaling $n\alpha_n \sim 1$ as $n\to \infty$, the term in square brackets
is of the same order as the Riemann sum for the integral $\int_0^1 V(t)dt < \infty$. Therefore, the condition \eqref{Riemann2} reduces to
\begin{equation}\label{previous}
\frac{K}{n\sigma h \alpha_n^2} \sim 1,
\end{equation}
which leads to the following expressions for $\alpha_n$ and $\ell_n$ (by using also \eqref{betan}):
\begin{equation}\label{betaCI}
\alpha_n^{(2)} \sim \sqrt{\frac{K}{n\sigma h}} = \beta_n, \qquad \ell_n^{(2)} \sim \sqrt{\frac{Knh}{\sigma}} = n\beta_n h.
\end{equation}
Note that by \eqref{betaCI} we can directly reformulate the scaling regime as $\beta_n\sim \frac{1}{n}$. More precisely, this corresponds to $\beta_n = \frac{c_n}{n}$, for some $c_n\sim 1$ ($c_n\to c$ as $n\to\infty$). The rescaling of the energy \eqref{Energy2} obtained in this case is therefore
\begin{align}\label{E2:d}
E_n^{(2)}(x):= \frac{c_n}{n^2}\sum_{k=1}^{n}\sum_{j=0}^{n-k} V\left(c_n(x_{j+k}-x_j)\right) + \frac1n\sum_{j=0}^n x_j.
\end{align}

\subsection{Continuum limit: derivation of the internal stress}
We rewrite the discrete energy \eqref{E2:d} in terms of the empirical measure $\rho_n= \frac{1}{n}\sum_{i=1}^n\delta_{x_i}$, and we obtain
\begin{equation}\label{E2:dn}
E_n^{(2)}(x)= \frac{c_n}{2}\int_0^\infty\int_0^\infty V(c_n(x-y))\rho_n(x)\rho_n(y)dx dy + \int_0^\infty x \rho_n(x)dx.
\end{equation}
For a large number of walls, i.e., as $n\to \infty$, $\rho_n$ converges to a continuum density $\rho$, $c_n\to c$, and the energy $E^{(2)}_n(x)$ converges to the continuum functional $E^{(2)}$ defined as:
\begin{equation}\label{E2:c}
E^{(2)}(\rho):= \frac c2\int_0^\infty\int_0^\infty V(c(x-y))\rho(x)\rho(y)dx dy + \int_0^\infty x \rho(x)dx.
\end{equation}

The equilibrium dislocation density $\rho$ is the solution of the Euler-Lagrange equation associated with the functional $E^{(2)}$, i.e., the solution of the integral equation
\begin{equation}\label{EL:2}
c\int_0^{\infty}V\big(c(x-y)\bigr)\partial_y \rho(y) dy + 1 = 0
\end{equation}
for every $x\in (0,\infty)$, or, equivalently,
\begin{equation}\label{EL:2a}
V_c\ast\,\partial_x\rho + 1 = 0,
\end{equation}
where $V_c(s) := cV(cs).$
This equation is again the mesoscopic equilibrium equation $\sigma^{(2)}_{\textrm{int}} - \sigma=0$ in its non-dimensional form.
Hence, the dimensionless internal stress obtained from this rescaling is
\begin{equation*}
\sigma^{(2)}_{\textrm{int}} = - V_c\ast\,\partial_x\rho.
\end{equation*}
The dimensional version, for $c=1$, is~\eqref{intstress2dimensional}.

We note that, as for the previous scaling regime, also in this case the continuum equilibrium equation is a singular integral equation. The main difference is that whereas in the rescaling case~(1) we could approximate $V$ with its limit behaviour near zero (namely its logarithmic behaviour), here the complete energy density $V$ enters the limit functional \eqref{E2:c} and the equilibrium equation \eqref{EL:2a}; the  scaling constant $c = \lim_{n\to\infty} n\beta_n$ enters the expression as well.

\subsection{Comparison Discrete vs Continuum}
\begin{figure}[htbp]
\labellist
\pinlabel {\small Dislocation wall positions} [b] at 180 -18
\pinlabel \small $\rho^d$ at 328 246.5
\pinlabel \small $\rho$ at 326 256
\endlabellist
\begin{center}
\includegraphics[width=3.6in]{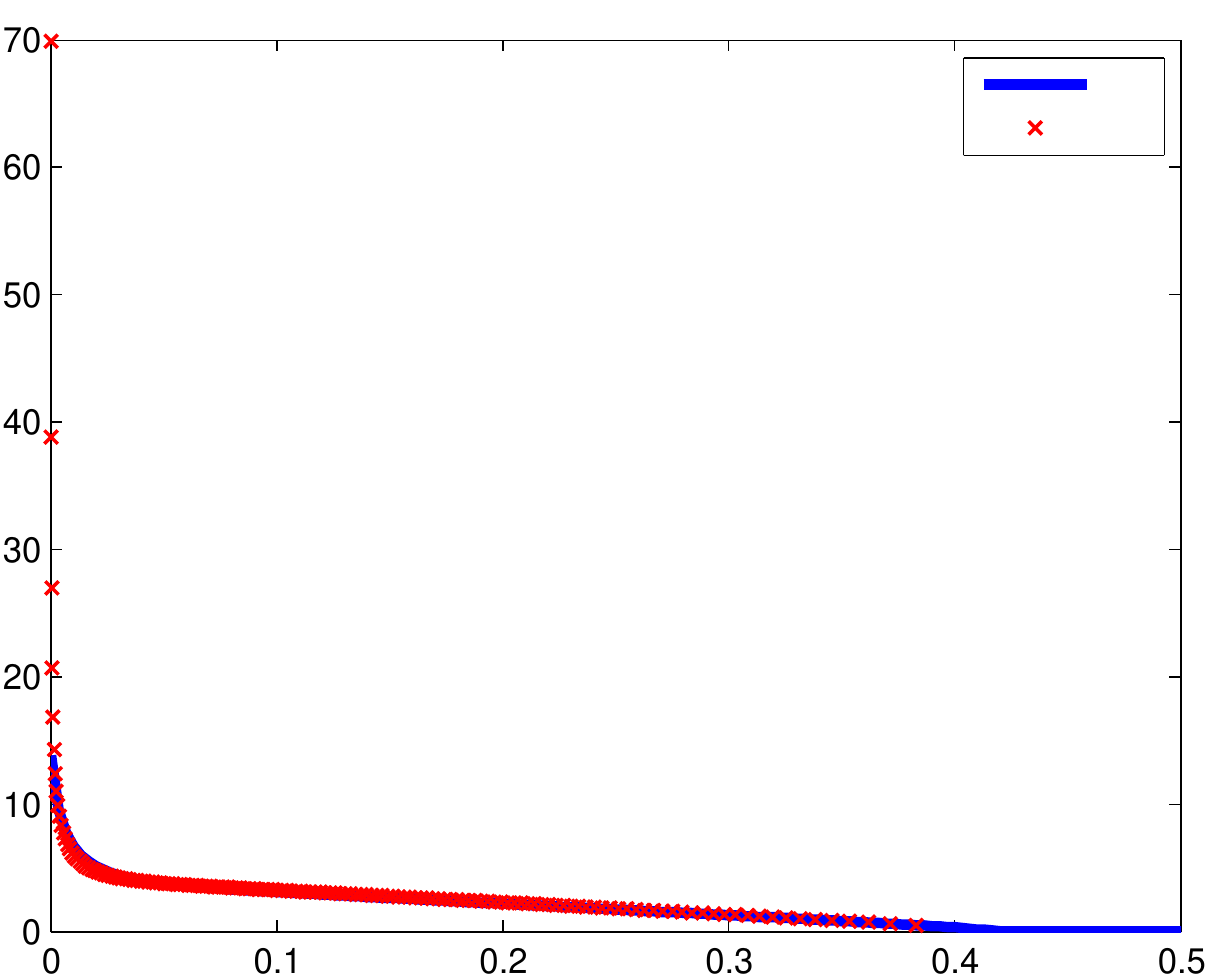}
\end{center}
\caption{Comparison of the discrete and continuous pile-ups for $E_n^{(2)}$ and $E^{(2)}$, for $n=150$ $\beta_n = 5/n$.}
\label{DD2}
\end{figure}

In this section we show the agreement between the solution of the upscaled continuum equation \eqref{EL:2} and the minimiser of the discrete energy \eqref{E2:d}, for a large number $n$ of dislocation walls. The plot is shown in Figure \ref{DD2}.

\section{Intermediate regime $\frac{1}{n}\ll \alpha_n \ll 1$: derivation of the upscaled internal stress}\label{Intermediate}

This scaling regime corresponds to the \textit{intermediate} situation in which the pile-up length $\ell_n$ is much larger than the vertical spacing $h$ while the average in-plane distance between consecutive dislocations is, on the contrary, much smaller than $h$.
This case corresponds to the \textit{remote region} considered in \cite{TRB}.

\subsection{Heuristics for the scaling of the discrete energy}
Using again the equispaced wall distribution $x_i=\frac{i}{n}$ as a test for the discrete energy $E_n$ in \eqref{Energy2} and proceeding as in the previous cases leads to the following requirement
\begin{align}\label{bound3}
1\sim\frac{K}{n^2\sigma h \alpha_n}\sum_{k=1}^{n}\sum_{j=0}^{n-k} V\left(\alpha_nk\right) \sim \frac{K}{n\sigma h \alpha_n^2} \left[\alpha_n\sum_{k=1}^{n} V\left(\alpha_nk\right)\right];
\end{align}
since, by assumption, $\alpha_n\to 0$ and $n\alpha_n\to \infty$, the term in the square brackets is a Riemann sum for the integral $\int_0^\infty V(t)dt =1/6\pi$. Hence the condition \eqref{bound3} on the energy is equivalent to the following relation among the parameters:
\begin{equation}\label{betaI}
\frac{K}{n\sigma h \alpha_n^2}\sim 1,
\end{equation}
which is identical to \eqref{previous} and therefore leads to the same expressions for $\alpha_n^{(3)}$ and $\ell_n^{(3)}$ as for $\alpha_n^{(2)}$ and $\ell_n^{(2)}$ in \eqref{betaCI}.
In particular, also in this case $\alpha_n^{(3)}\sim \beta_n$; in terms of $\beta_n$ the scaling regime is $\frac{1}{n}\ll\beta_n\ll1$.
This relation between $\ell_n^{(3)}$ and the parameters $h$, $\sigma$, $K$ and $n$ found above was also found in \cite{TRB}.

\medskip

Using \eqref{betaI}, the rescaling of the energy \eqref{Energy2} becomes
\begin{align}\label{E3:d}
E_n^{(3)}(x):= \frac{\beta_n}{n}\sum_{k=1}^{n}\sum_{j=0}^{n-k} V\left(n\beta_n(x_{j+k}-x_j)\right) + \frac1n\sum_{j=0}^n x_j,
\end{align}
for $\frac{1}{n}\ll\beta_n\ll1$.

\subsection{Continuum limit: derivation of the internal stress}
In terms of the distribution $\rho_n= \frac{1}{n}\sum_{i=1}^n\delta_{x_i}$, the energy \eqref{E3:d} can be rewritten as
\begin{equation}\label{E3:dn}
E_n^{(3)}(x)= \frac12\,n\beta_n\int_0^\infty\int_0^\infty V\left(n\beta_n(x-y)\right)\rho_n(x)\rho_n(y)dx dy + \int_0^\infty x \rho_n(x)dx.
\end{equation}
We note that, since $n\beta_n\to \infty$, $n\beta_nV(n\beta_n s)\to \left(\int_{-\infty}^\infty V\right)\delta_0 = (3\pi)^{-1}\delta_0$; hence for $n\to\infty$ the energy $E^{(3)}_n$ converges to the continuum energy $E^{(3)}$ defined as:
\begin{equation}\label{E3:c}
E^{(3)}(\rho):=\frac1{6\pi}\int_0^\infty \rho^2(x)dx + \int_0^\infty x\rho(x)dx.
\end{equation}
The dislocation density $\rho$ minimising $E^{(3)}$ is the solution of the Euler-Lagrange equation associated to the dimensionless functional $E^{(3)}$, i.e.,
\begin{equation}\label{EL:3}
\frac1{3\pi}\partial_x\rho + 1 = 0.
\end{equation}
We notice that unlike the previous scaling regimes the continuum equilibrium equation in this case is local.
Moreover, from \eqref{EL:3} it is clear that the optimal dislocation density is linear. This was observed numerically in~\cite{TRB} and proved in~\cite{Hall11} using formal methods.

\medskip

Equation \eqref{EL:3} is the mesoscopic equilibrium equation $\sigma^{(3)}_{\textrm{int}} - \sigma=0$ in its non-dimensional form.
Hence, the dimensionless internal stress obtained from this rescaling is
\begin{equation*}
\sigma_{\textrm{int}}^{(3)} = -\frac1{3\pi}\partial_x \rho.
\end{equation*}

\subsection{Comparison Discrete vs Continuum}
\begin{figure}[htbp]
\labellist
\pinlabel {\small Dislocation wall positions} [b] at 180 -18
\pinlabel \small $\rho^d$ at 328 246.5
\pinlabel \small $\rho$ at 326 256
\endlabellist
\begin{center}
\includegraphics[width=3.6in]{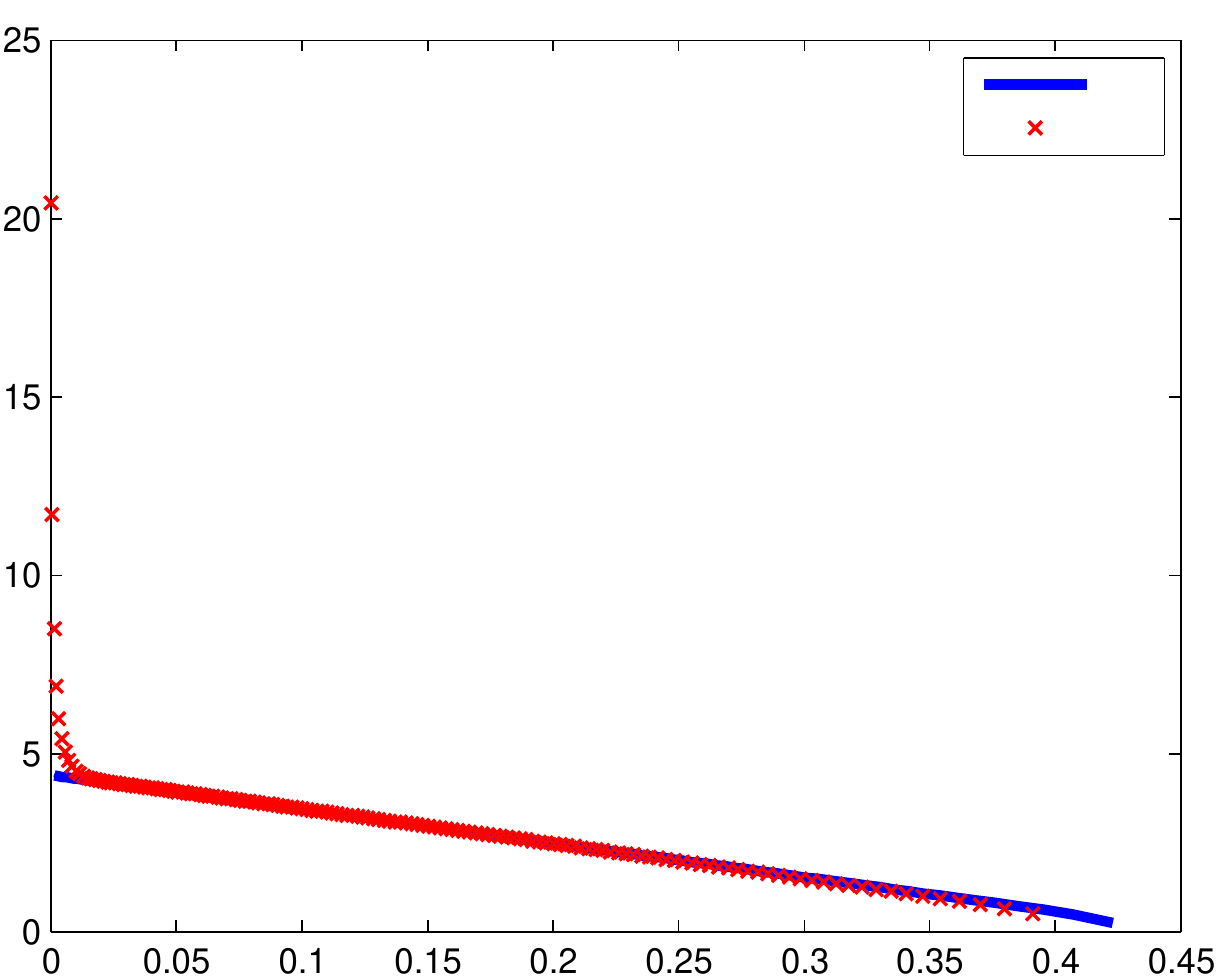}
\end{center}
\caption{Comparison of the discrete and continuous pile-ups for $E_n^{(3)}$ and $E^{(3)}$, for $n=150$ and $\beta_n = 1/\sqrt n=1/\sqrt{150}$.}
\label{Fig:linear}
\end{figure}

In Figure \ref{Fig:linear} we compare numerically the discrete density corresponding to a minimiser of the rescaled energy $E_n^{(3)}$ in \eqref{E3:d} and the solution of the continuum equation \eqref{EL:3}. We notice that the agreement between the two densities is perfect in the \textit{bulk}, sufficiently far from the boundary of the pile-up region, where the occurrence of boundary layers is expected.

\subsection{Comparison with the Evers-Brekelmans-Geers continuum model.}

The internal stress $\sigma^{(3)}_{\textrm{int}}$ coincides, at least qualitatively, with the one proposed in~\cite{EvGeers} and derived phenomenologically from interactions among dislocations. In the quoted paper, however, the expression of the internal stress contains a length scale $R$ (not explicitly determined), representing the spatial reach of the dislocations interactions. By matching the stress $\sigma^{(3)}_{\textrm{int}}$ in its dimensional form \eqref{intstress3dimensional} with the one derived in \cite{EvGeers} we can determine the length scale $R$, namely $R\sim h$.
Therefore $R$ depends on the in-wall spacing $h$ only.

\section{Second critical regime $\alpha_n\sim 1$: derivation of the upscaled internal stress}\label{SecondCritical}
This scaling regime corresponds to configurations where the in-plane and the in-wall distances between consecutive dislocations are of the same order. As a result, the length of the pile-up region is much larger than $h$.

The scaling $\alpha_n\sim 1$ is a special situation as compared to the cases (1)--(3) treated so far and the last case (5). In fact in all the cases (1)--(3) the aspect ratio is very small ($\alpha_n\ll1$), and so the in-plane interaction is stronger than the in-wall interaction, while the opposite situation arises in case (5), since $\alpha_n\gg1$. For this reason we look at case (4) as the main critical case and we expect to get a better insight in the problem from its analysis.

\subsection{Heuristics for the scaling of the discrete energy}
Proceeding exactly as in the previous cases leads to the same requirement as in \eqref{bound3}, where this time $\alpha_n\sim 1$.
Rewriting the condition~\eqref{bound3} for an equidistant distribution of dislocations we have
\begin{equation*}
1\sim\frac{K}{n^2 \sigma h \alpha_n}\sum_{k=1}^{n}\sum_{j=0}^{n-k} V(\alpha_n k) \sim \frac{K}{n \sigma h \alpha_n}\sum_{k=1}^{n}V(\alpha_n k).
\end{equation*}
Hence, since $\sum_{k=1}^{\infty}V(k)<\infty$, the boundedness of the energy reduces to the following condition on the coefficients:
\begin{equation}\label{betaIR}
\frac{K}{n\sigma h \alpha_n} \sim 1.
\end{equation}
From this bound we obtain the following expressions for the aspect ratio $\alpha_n$ and the pile-up length~$\ell_n$:
\begin{equation}\label{betaIR2}
\alpha_n^{(4)} = \frac{K}{n\sigma h} = \beta_n^2, \quad \ell_n^{(4)} = \frac{K}{\sigma}.
\end{equation}
We note that the scaling regime can be rephrased as $\beta_n \sim 1$. Therefore
we can rescale the energy \eqref{Energy2} as in \eqref{E3:d}, i.e.,
\begin{equation}\label{E4:d}
E_n^{(4)}(x):= \frac{\beta_n}{n}\sum_{k=1}^{n}\sum_{j=0}^{n-k} V\left(n\beta_n(x_{j+k}-x_j)\right) + \frac1n\sum_{j=0}^n x_j,
\end{equation}
where now $\beta_n\sim 1$ (note that \eqref{betaIR2} suggests $\beta_n^2$ instead of $\beta_n$, but in this regime $\beta_n\sim \beta_n^2$).

\subsection{Continuum limit: derivation of the internal stress}
The derivation of the continuum energy in this regime is quite different from the one outlined in Sections \ref{Subcritical}-\ref{Intermediate}. This can be explained by considering the discrete energy \eqref{Energy2}. While for $\alpha_n\ll1$ (namely in the scaling regimes (1)--(3)) the sum over $k$ gives rise to an integral term, in addition to the integral generated by the sum over~$j$, for $\alpha_n\sim 1$ this is no longer the case, so a different approach has to be followed.
We mention here the idea only briefly (setting $\beta_n=1$ for simplicity) and refer to \cite{GPPS_Math} for the complete proof.

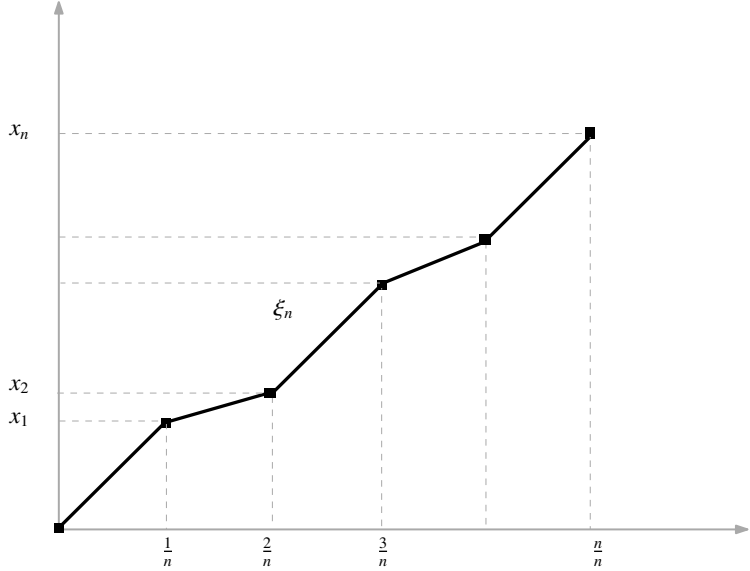
\begin{figure}[htbp]
\labellist
\pinlabel {\small $\frac 1n$} [b] at 160 100
\pinlabel {\small $\frac 2n$} [b] at 220 100
\pinlabel {\small $\frac 3n$} [b] at 290 100
\pinlabel {\small $\frac nn$} [b] at 420 100
\pinlabel {\small $x_1$} [b] at 70 185
\pinlabel {\small $x_2$} [b] at 70 205
\pinlabel {\small $x_n$} [b] at 70 360
\pinlabel {\small $\xi_n$} [b] at 230 250
\endlabellist
\begin{center}
\includegraphics[width=3.6in]{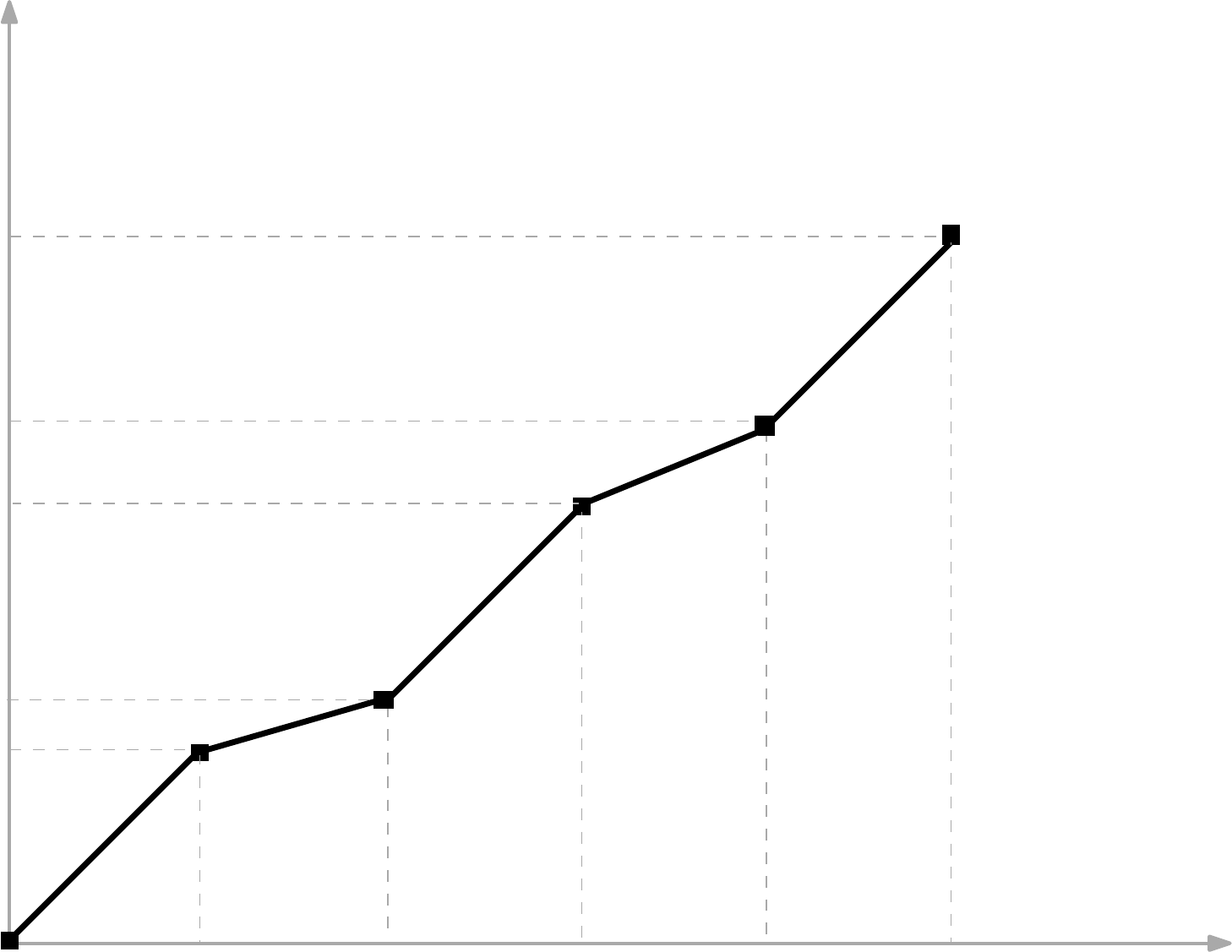}
\end{center}
\caption{Interpolation of the positions of the walls.}
\label{Fig:interp}
\end{figure}

The idea in this case is to view the positions of the walls $x_i$ as the deformed positions from an initial equispaced wall configuration of $n+1$ walls in $(0,1)$ via a deformation map $\xi_n$, i.e., $x_i=\xi_n\left(\frac{i}{n}\right)$. The map $\xi_n$ can be extended in the whole range $(0,1)$ as the continuum affine interpolation of $x_1,\dots,x_n$ (see Figure \ref{Fig:interp}).

In this way we can formally read the argument of $V$ in \eqref{E4:d} as a quotient ratio of $\xi_n$, i.e.,
$$
 V\left(n(x_{j+k}-x_j)\right) =  V\left(k\,\frac{x_{j+k}-x_j}{k/n}\right) \sim V\left(k\,\xi'_n\left(\frac{j}{n}\right)\right).
$$
Hence for the discrete energy \eqref{E4:d} we have
\begin{equation}\label{E4rewrite}
E_n^{(4)}(x)= \sum_{k=1}^{n}\left\{\frac{1}{n}\sum_{j=0}^{n-k} V\left(k\,\xi'_n\left(\frac{j}{n}\right)\right)\right\} + \frac1n\sum_{j=0}^n \xi_n\left(\frac{j}{n}\right),
\end{equation}
which for $n\to \infty$ converges to the continuum functional
\begin{equation}\label{E4:x}
E^{(4)}(\xi):= \sum_{k=1}^{\infty}\int_0^1 V(k \xi'(s))\,ds + \int_0^1\xi(s)\,ds,
\end{equation}
where $\xi$ is the limit of the interpolations and the integrals are the limits of the Riemann sums in \eqref{E4rewrite}. Since we are interested in a continuum model in terms of the dislocation density rather than the limit positions it remains to change variables in the energy \eqref{E4:x}, using the relation
\begin{equation}\label{cdensity}
\rho(\xi(s)):= \frac{1}{\xi'(s)},
\end{equation}
which is the continuum version of the discrete relation \eqref{ddensity}. The relation \eqref{cdensity} entails $ds = \rho(x)\,dx$. Changing variables in \eqref{E4:x} leads to a continuum energy depending on the dislocation density, namely
\begin{equation}\label{E4:c}
E^{(4)}(\rho)=\int_{0}^\infty V_{\textrm{eff}}\left(\frac{1}{\rho(x)}\right)\rho(x) dx + \int_0^\infty x\rho(x)dx,
\end{equation}
where $V_{\textrm{eff}}(t):= \sum_{k=1}^{\infty}V(kt)$, for every $t\in \mathbb{R}$.

The dislocation density $\rho$ minimising the energy \eqref{E4:c} is the solution of the Euler-Lagrange equation associated to the dimensionless functional $E^{(4)}$, i.e.,
\begin{equation}\label{EL:4}
-\frac{\partial_x\rho}{\rho^3}\,\varphi'_{\textrm{eff}}\left(\frac{1}{\rho}\right) + 1 = 0,
\end{equation}
where $\varphi_{\textrm{eff}}(t):= \sum_{k=1}^{\infty}k\varphi(kt)$, for every $t\in \mathbb{R}$.
Note that the infinite sum that appears in the definition of the continuum equilibrium equation \eqref{EL:4} results form having taken into account all the dislocation interactions in the discrete model. More precisely, the term $k=1$ in the sum corresponds to the contribution of the interactions among nearest neighbours, the term $k=2$ to next-to-nearest neighbour interactions and so on. Moreover, whereas in the previous cases the dislocation walls were sufficiently close to regard them as a continuum density, here the discrete interactions prevail.

Equation \eqref{EL:4} is again the mesoscopic equilibrium equation $\sigma^{(4)}_{\textrm{int}} - \sigma=0$ in its non-dimensional form.
Hence, the dimensionless internal stress obtained for this rescaling is
\begin{equation*}
\sigma^{(4)}_{\textrm{int}} = \frac{\partial_x \rho}{\rho^3}\,\varphi'_{\textrm{eff}}\left(\frac{1}{\rho}\right).
\end{equation*}
The dimensional version is~\eqref{intstress4dimensional}.

\begin{rem}
We show here how \eqref{EL:4} can be heuristically derived starting from the original, dimensional equilibrium equations \eqref{Eq:1}. We will focus on a general $k$-th term in the definition of~$\varphi_{\textrm{eff}}$.

For a dislocation located at a point $\tilde x$, we denote by $d_k$ the average distance to its $k$-th neighbours, i.e., $d_k\sim \frac{k}{\tilde \rho}$. The actual distance $d_k^{\pm}(\tilde x)$ between the dislocation at $\tilde x$ and its $k$-th neighbour on the right and on the left are given by the following corrections of $d_k$, namely
\begin{align*}
d_k^+ &\sim \frac{k}{\tilde \rho} + \left(-\frac{1}{\tilde \rho^2}\frac{\partial\tilde \rho}{\partial \tilde x}\right)\left(\frac12 \,\frac{k}{\tilde \rho}\right)\\
d_k^- &\sim \frac{k}{\tilde \rho} + \left(-\frac{1}{\tilde\rho^2}\frac{\partial\tilde \rho}{\partial \tilde x}\right)\left(-\frac12 \,\frac{k}{\tilde \rho}\right).
\end{align*}

\begin{figure}[htbp]
\labellist
\pinlabel {\small $\tilde x$} at 290 235
\pinlabel {\small $d_k^-$} at 200 280
\pinlabel {\small $d_k^+$} at 360 280
\endlabellist
\begin{center}
\includegraphics[width=3in]{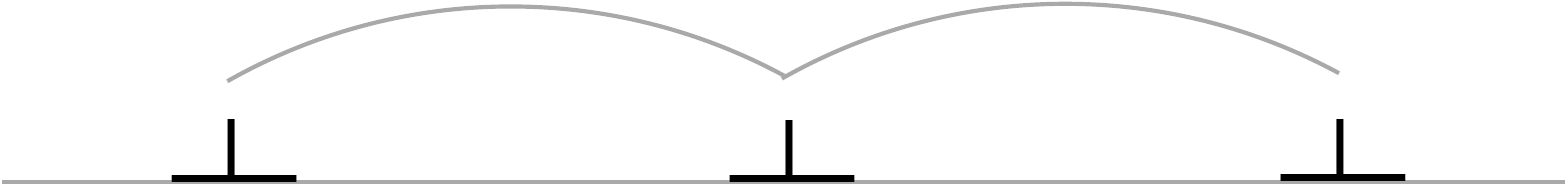}
\end{center}
\caption{$k$-interactions.}
\label{Ron}
\end{figure}

Therefore, the distance between its $k$-th neighbours is
\[
\Delta d \sim -\frac{k}{\tilde \rho^3} \partial_{\tilde x} \tilde \rho.
\]
According to the previous relations, the force exerted on the dislocation at $\tilde x$ by the $k$-th neighbour dislocation on the right and on the left are, respectively,
\begin{align*}
\frac Kh \varphi\Big(\frac{d_k^+}h\Big) &\sim \frac Kh\bigg[\varphi\Big(\frac{d_k}h\Big) + \frac1{2h}\,\varphi'\Big(\frac{d_k}h\Big)\Delta d\bigg],\\
\frac Kh \varphi\Big(\frac{d_k^-}h\Big) &\sim \frac Kh\bigg[\varphi\Big(\frac{d_k}h\Big) - \frac1{2h}\,\varphi'\Big(\frac{d_k}h\Big)\Delta d\bigg].
\end{align*}
The net force exerted on the dislocation at $\tilde x$ by its $k$-th neighbours is, therefore,
\[
\frac Kh\Delta \varphi \sim  \frac K{h^2} \varphi'(d_k)\Delta d = - \frac K{h^2}\frac{k}{\tilde \rho^3}\, \partial_{\tilde x} \tilde \rho\;\varphi'\left(\frac{k}{\tilde \rho}\right),
\]
which is exactly the $k$-th term in \eqref{intstress4dimensional}.
\end{rem}

\subsection{Comparison Discrete vs Continuum}
\begin{figure}[htbp]
\labellist
\pinlabel {\small Dislocation wall positions} [b] at 180 -17
\pinlabel \small $\rho^d$ at 322 246.5
\pinlabel \small $\rho$ at 318 255
\endlabellist
\begin{center}
\includegraphics[width=3.6in]{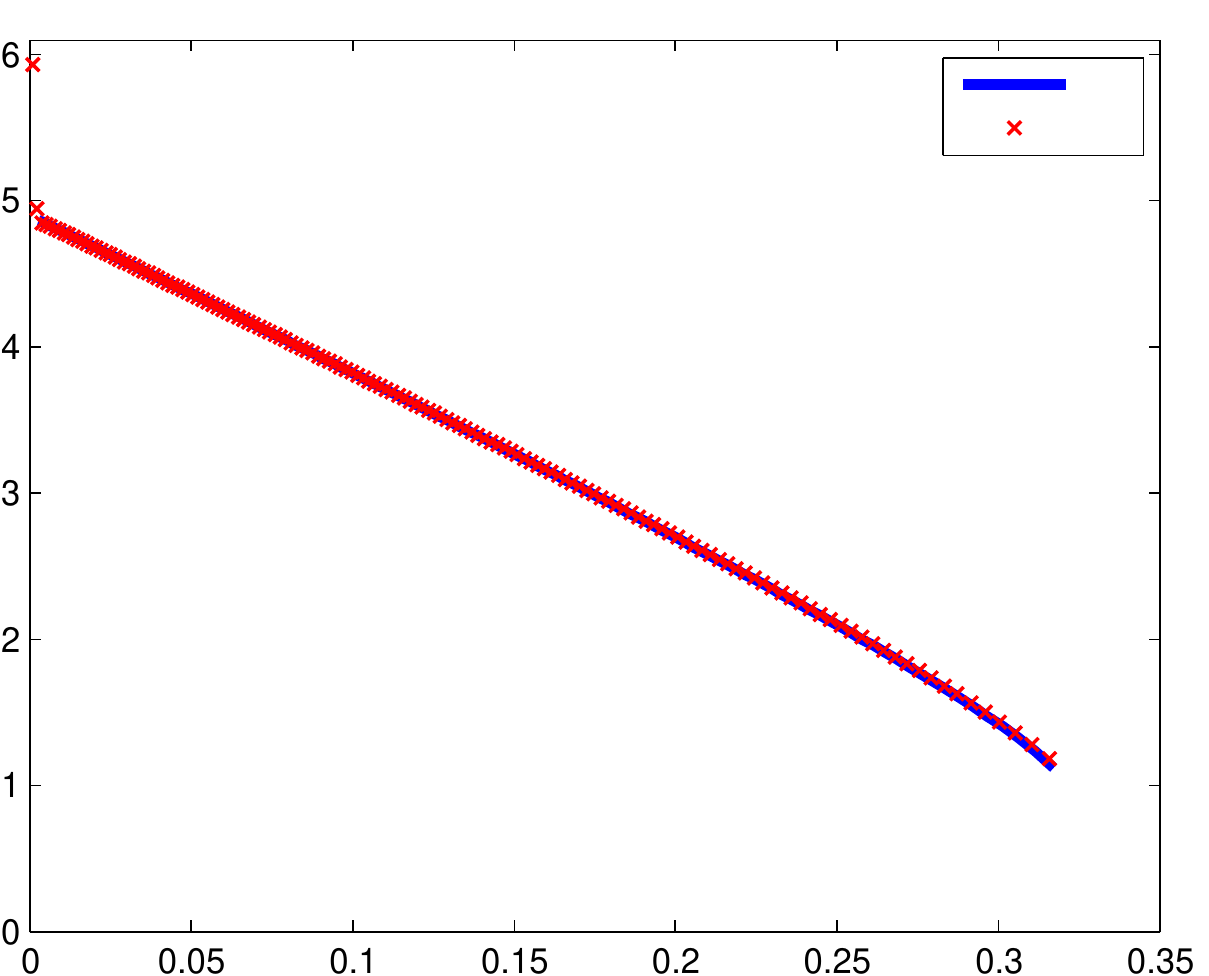}
\end{center}
\caption{Comparison of the discrete and continuous pile-ups for $E_n^{(4)}$ and $E^{(4)}$, with $n=150$ and $\beta_n =1$.}
\label{Case4}
\end{figure}

In this section we show the agreement between the solution of the continuum equation \eqref{EL:4} and the minimiser of the rescaled discrete energy \eqref{E4:d}, for $\beta_n\sim 1$. The agreement is shown in Figure \ref{Case4}; also in this case we can notice the occurrence of boundary layers at the left end of the pile-up region.

\section{Supercritical regime $\alpha_n\gg 1$: derivation of the upscaled internal stress}\label{Supercritical}
In this scaling regime the in-plane dislocation distance is much larger than the in-wall distance, and therefore the in-plane interaction is quite weak.
In this case we can prove that truncating the in-plane interactions to the first neighbours (but keeping the walls) leads to the right
limit model, unlike the previous cases, where all interactions had to be accounted for.

\subsection{Heuristics for the scaling of the discrete energy}
We first note that the energy \eqref{Energy2} is always larger than its truncation to first neighbours, i.e.,
\begin{equation}\label{Trunc1}
E_n(x)\geq\frac{K}{n^2 \sigma h \alpha_n}\sum_{j=0}^{n-1} V\left(n\alpha_n(x_{j+1}-x_j)\right) + \frac1n\sum_{j=0}^n x_j,
\end{equation}
regardless of the scaling regime. In the case $\alpha_n \gg 1$  the bound \eqref{Trunc1} is optimal, meaning that the energy $E_n$ and its first order truncation give rise to the same continuum model (see \cite{GPPS_Math}). This is in contrast with what happens in \textit{all} the previous cases (1)--(4), where the truncation of the interactions would produce a completely different (and wrong) result (see also the discussion in~\cite{RPG}).

To obtain the energy scaling we use once more the equispaced wall configuration $x_i=\frac{i}{n}$ as a test configuration. Then the bound on the energy (actually the bound on the energy truncated to the first neighbours) reduces to

\begin{equation}\label{C:1}
\frac{K}{n \sigma h \alpha_n}\, V(\alpha_n)\sim 1.
\end{equation}
Since $\alpha_n \gg 1$, we can substitute for $V$ its asymptotic behaviour at infinity; by \eqref{defV} we have
\begin{equation}\label{V:inf}
V(s) \sim \frac{2}{\pi}s e^{-2\pi s} \quad \mbox{as}\quad s\to \infty.
\end{equation}
Therefore, using the previous relation, the bound \eqref{C:1} becomes
\begin{equation*}
\frac{K}{n \sigma h}\frac{2}{\pi} e^{-2\pi\alpha_n} \sim 1.
\end{equation*}
In terms of $\alpha_n$, and using \eqref{ab1}, the previous relation entails
\begin{equation}\label{alb}
\alpha_n^{(5)} \sim \frac{1}{2\pi}\,\log\left(\frac{2}{\pi}\frac{K}{nh\sigma}\right) = \frac{1}{2\pi}\log\left(\frac{2 \beta_n^2}{\pi}\right), \qquad
\ell_n^{(5)} = nh\alpha_n \sim \frac{nh}{2\pi}\log\left(\frac{2 \beta_n^2}{\pi}\right).
\end{equation}
Therefore, the energy scaling for $\alpha_n\gg1$ can be rephrased in terms of $\beta_n$ as $\beta_n\gg1$, and the corresponding energy scaling in the regime $\beta_n\gg1$ is
\begin{equation}\label{E5:d}
E_n^{(5)}(x):= \frac{2\pi\beta_n^2}{n\log\left(\frac{2\beta_n^2}{\pi}\right)}
\sum_{k=1}^{n}\sum_{j=0}^{n-k} V\left(\frac{n}{2\pi}\,\log \left(\frac{2 \beta_n^2}{\pi}\right) (x_{j+k}-x_j)\right) + \frac1n\sum_{j=0}^n x_j.
\end{equation}
For what follows it is convenient to introduce a new parameter $\gamma_n:= \frac{1}{2\pi}\,\log \left(\frac{2 \beta_n^2}{\pi}\right)$ (note that this scaling regime corresponds to $\gamma_n\gg1$) and to rewrite \eqref{E5:d} in terms of $\gamma_n$ as follows:
\begin{equation}\label{E5:dg}
E_n^{(5)}(x)= \frac{\pi}{2}\frac{e^{2\pi\gamma_n}}{n\gamma_n}
\sum_{k=1}^{n}\sum_{j=0}^{n-k} V\left(n\gamma_n(x_{j+k}-x_j)\right) + \frac1n\sum_{j=0}^n x_j.
\end{equation}

\subsection{Continuum limit: derivation of the internal stress}
As in Section \ref{SecondCritical} we introduce the interpolation $\xi_n$ of the positions of the walls. Heuristically we expect the term $k=1$ in the sum in \eqref{E5:dg} to be dominant, and hence we disregard the other terms in our formal derivation of the limit energy (for the rigorous proof we refer to \cite{GPPS_Math}). We have, formally,
\begin{align*}
E_n^{(5)}(x)&\simeq \frac{\pi}{2}\frac{e^{2\pi\gamma_n}}{n\gamma_n}\sum_{j=0}^{n-1} V\left(\gamma_n\xi'\left(\frac{j}{n}\right)\right) + \frac1n\sum_{j=0}^n \xi\left(\frac{j}{n}\right)\nonumber\\
& \simeq \frac{\pi}{2}\frac{e^{2\pi\gamma_n}}{\gamma_n}\int_0^1 V(\gamma_n\xi'(s))ds + \int_0^1 \xi(s)ds.
\end{align*}
By \eqref{V:inf} we have that, for $s\in (0,1)$,
$$
\frac{\pi}{2}\frac{e^{2\pi\gamma_n}}{\gamma_n} \, V(\gamma_n\xi'(s))\sim \xi'(s)\, e^{2\pi(\gamma_n-\xi'(s))},
$$
which is finite (actually zero) only if $\xi'(s)\geq 1$. In terms of the density $\rho$, using \eqref{cdensity} the limit energy is the dimensionless continuum functional $E^{(5)}$ defined as:
\begin{equation}\label{E5:c}
E^{(5)}(\rho) = \begin{cases}
\smallskip
\displaystyle \int_0^\infty x \rho(x)dx \quad &\mbox{if } \rho\leq 1, \\
+\infty \quad &\mbox{otherwise}.
\end{cases}
\end{equation}

This limit functional is degenerate: it only is finite if the dislocations are sufficiently far apart. The optimal density in this case is the constant density $1$ from $x=0$ to $x=1$, and density zero for $x\geq 1$ (see Figure~\ref{Fig:constant} below).

\subsection{Comparison Discrete vs Continuum}
In this section we show the agreement between the minimiser of the discrete energy \eqref{E5:d}, for $\beta_n \gg 1$ and for large $n$, and the minimiser of the continuum energy \eqref{E5:c}.

\begin{figure}[htbp]
\labellist
\pinlabel {\small Dislocation wall positions} [b] at 180 -17
\pinlabel \small $\rho^d$ at 328 243.5
\pinlabel \small $\rho$ at 326 253
\endlabellist
\begin{center}
\includegraphics[width=3.6in]{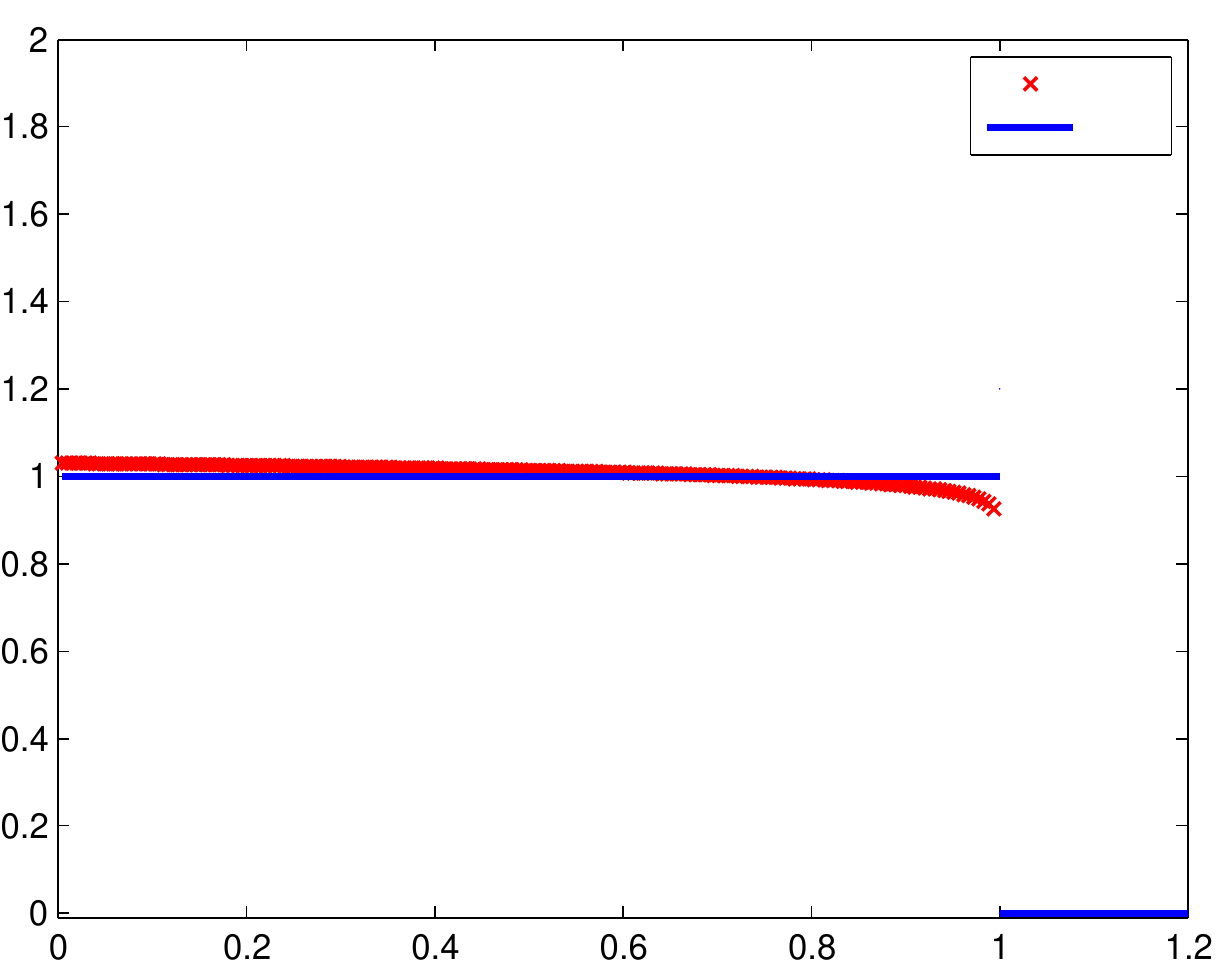}
\end{center}
\caption{Optimal densities for $E_n^{(5)}$ and $E^{(5)}$, with $n=150$ and $\beta_n =10^8$.}
\label{Fig:constant}
\end{figure}

In Figure \ref{Fig:constant} we have plotted the piecewise constant optimal continuum density and the optimal discrete density for $n=150$ and $\beta_n = 10^8$: we see that the discrete density approaches the value~$1$.

\section{Comparison with previous models}\label{Comparison}

In this section we show how the expressions for the internal stress obtained in this paper, starting from an idealised discrete model, relate to some well-known models in the engineering literature. Our results offer therefore a unifying approach for understanding these existing models: They can all be derived by upscaling the same discrete model, but under different assumptions on the local arrangements of the dislocations, i.e., for different values of the aspect ratio $\beta_n$.

\medskip

First of all we list the dimensional internal stresses obtained in the scaling regimes (1)--(4) (we omit (5) since the equilibrium equation is too degenerate to provide a useful expression for the internal stress), i.e.,

\begin{align*}
&(1)\,\, \mbox{Subcritical regime }& &\beta_n\ll \frac{1}{n}: & &\sigma^{(1)}_{\textrm{int}}(\tilde x)  =  \frac{K}{\pi^2}\int_0^\infty \log|\tilde x - \tilde y|\partial_{\tilde y}\tilde \rho(\tilde y)\, d\tilde y;\\
\bigskip
&(2)\,\, \mbox{First critical regime }& &\beta_n\sim \frac{1}{n}: & &\sigma^{(2)}_{\textrm{int}}(\tilde x) = - K \int V\left(\frac{\tilde x-\tilde y}{h}\right)\partial_{\tilde y}\tilde \rho(\tilde y)\,d\tilde y;
\end{align*}
\begin{align*}
&(3)\,\, \mbox{Intermediate regime }& &\frac{1}{n} \ll\beta_n\ll 1: & &\sigma^{(3)}_{\textrm{int}}(\tilde x) = -\frac{K h}{3\pi}\partial_{\tilde x}\tilde \rho;\\
\bigskip
&(4)\,\, \mbox{Second critical regime}& &\beta_n \sim 1: & &\sigma^{(4)}_{\textrm{int}}(\tilde x) = \frac{K}{h^2}\frac{1}{\tilde \rho^3}\, \varphi'_{\textrm{eff}}\left(\frac{1}{h\tilde\rho}\right)\, \partial_{\tilde x}\tilde \rho.
\end{align*}

\medskip

The back-stress $\sigma^{(1)}_{\textrm{int}}$ coincides with the one derived by Eshelby, Frank, and Nabarro~\cite{EFN} and Head and Louat~\cite{HL}. The starting point of both works is a system of discrete equilibrium equations for~$n$ dislocations in one slip plane, which is different from the \emph{walls} of dislocations that we consider. This translates into a system of discrete equilibrium equations of the form \eqref{Eq:1}, with the single-slip-plane interaction potential $\psi(t) = 1/(\pi^2t)$ instead of the the potential $\varphi$. We note that the interaction potential $\psi$ is similar to $\varphi$ for small values of $t$, but differs for larger values.

This is consistent with the fact that the internal stress $\sigma^{(1)}_{\textrm{int}}$ has been obtained in the scaling regime $\beta_n\ll 1/n$, where all the dislocation walls are confined in a region that is small compared to the slip planes spacing $h$. Hence the in-plane interactions between dislocations are much stronger than the in-wall interactions. An approximation by a single-slip-plane setup is therefore appropriate. A comparison between the two discrete models in terms of the equilibrium densities is shown in Figure \ref{EFNvsOurs}.

\medskip

The back-stress $\sigma^{(3)}_{\textrm{int}}$ coincides, at least qualitatively, with the one proposed by Evers, Brekelmans and Geers in \cite{EvGeers}.
In the quoted paper the authors derive phenomenologically an expression for the internal stress that is, like ours, the gradient of the dislocation density, up to a multiplicative constant (representing, in their case, the range of interaction of a dislocation).

In the special case $\beta_n = cn^{-1/2}$ for some constant $c$, the simple ordinary differential equation \eqref{EL:3} characterising the equilibrium dislocation density in this regime has also been obtained (following a different method) in the recent paper \cite{Hall11}, starting from the discrete system \eqref{Eq:1}.

\medskip

The other internal stresses $\sigma_{\mathrm{int}}^{(2)}, \sigma_{\mathrm{int}}^{(4)}$ and $\sigma_{\mathrm{int}}^{(5)}$ were not obtained so far.

\medskip

The expression for the internal stress proposed by Groma, Csikor, and Zaiser~\cite{Groma} is a special case. It has been derived under the assumption that the distance to a nearest neighbouring dislocation is independent of the direction in which it is found. In our formulation this corresponds to the second critical case (4), namely $\beta_n\sim 1$. Therefore it is interesting to compare the internal stress of~\cite{Groma} to $\sigma_{\mathrm{int}}^{(4)}$ above.

As it turns out, the internal stress of~\cite{Groma}, which up to constants reads $\sigma^{(GCZ)}_{\textrm{int}} = -\partial_{\tilde x}\tilde\rho/\tilde\rho$, can be \emph{formally} obtained from $\sigma_{\textrm{int}}^{(4)}$ by making two approximations.
The first  consists in truncating the number of interacting dislocations to the nearest neighbours.
As described above, in the discrete model of this paper we take into account all interactions. This can also be recognized in $\sigma_{\textrm{int}}^{(4)}$, which is defined in terms of $\varphi_{\mathrm{eff}}$, which in turn is a sum whose $k$-th term corresponds to the interactions
between $k$-th neighbours. Hence, truncating the interactions to the first neighbours is equivalent to replacing the effective potential
$\varphi_{\textrm{eff}}$ with $\varphi$. This first simplification reduces the internal stress $\sigma_{\textrm{int}}^{(4)}$ above to

\begin{equation}\label{step1}
\sigma_{\textrm{int}} = \frac{K}{h^2}\frac{1}{\tilde \rho^3}\, \varphi'\left(\frac{1}{h\tilde\rho}\right)\, \partial_{\tilde x}\tilde \rho.
\end{equation}

The second approximation that leads to the internal stress proposed in~\cite{Groma} is to substitute the force $\varphi(t)$ with its first-order Taylor-Laurent expansion close to zero, which coincides with the single-slip plane force $\psi(t) = \frac{1}{\pi^2 t}$ used in \cite{EFN,HL}.
By using $\psi$ instead of $\varphi$ in \eqref{step1} we obtain exactly (up to a constant)  $\sigma^{(GCZ)}_{\textrm{int}} = -\partial_{\tilde x}\tilde\rho/\tilde\rho$.

\begin{figure}[htbp]
\labellist
\pinlabel \scriptsize $\varphi$ at 333 246
\pinlabel \scriptsize $\varphi_{\textrm{eff}}$ at 335 259
\pinlabel \scriptsize $\psi$ at 333 235
\pinlabel $\ln s$ at 180 -5
\endlabellist
\begin{center}
\includegraphics[width=3.6in]{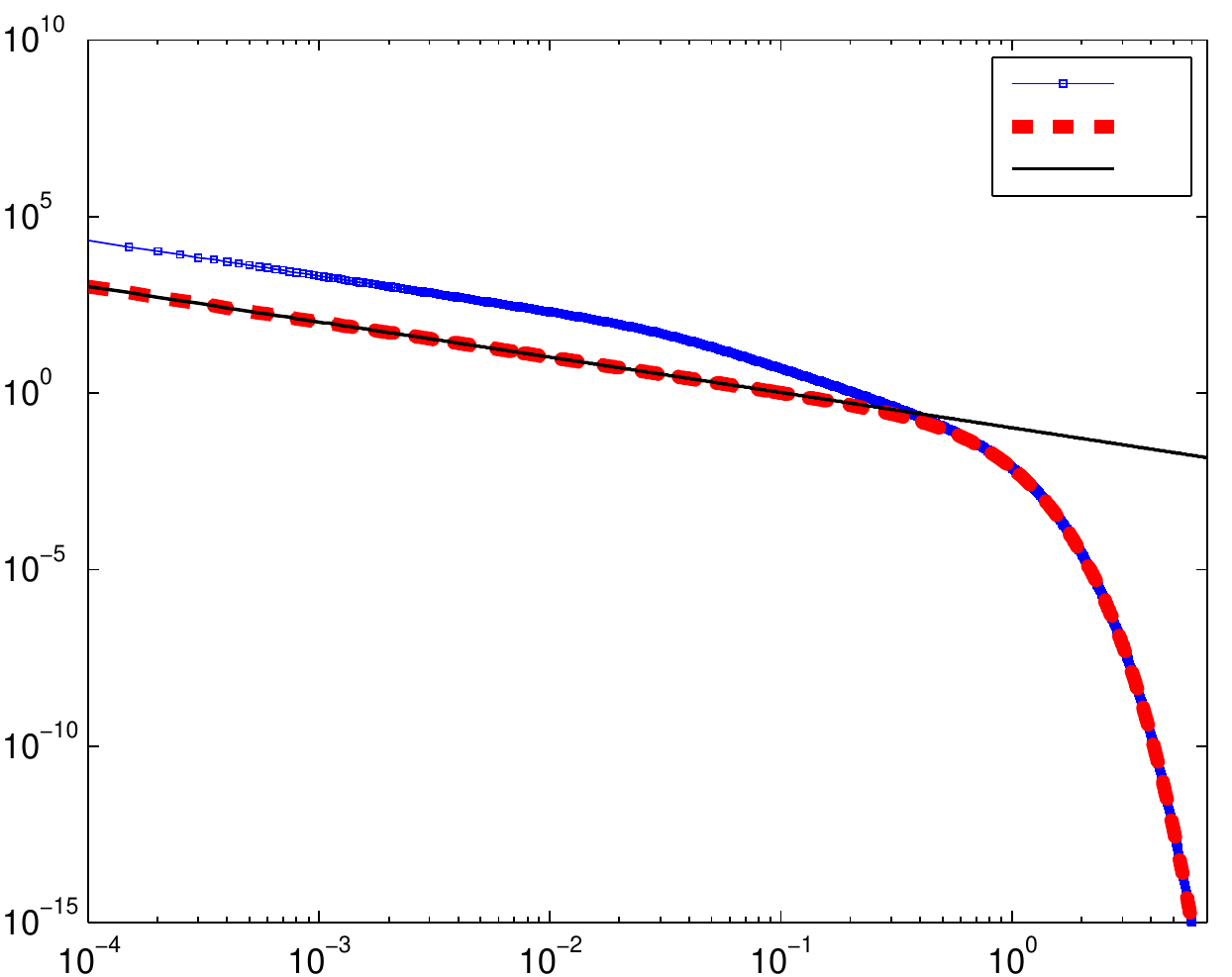}
\end{center}
\caption{Interaction potentials $\varphi_{\textrm{eff}}$, $\varphi$ and $\psi$ in a log-log plot.}
\label{Fig:Interactions}
\end{figure}
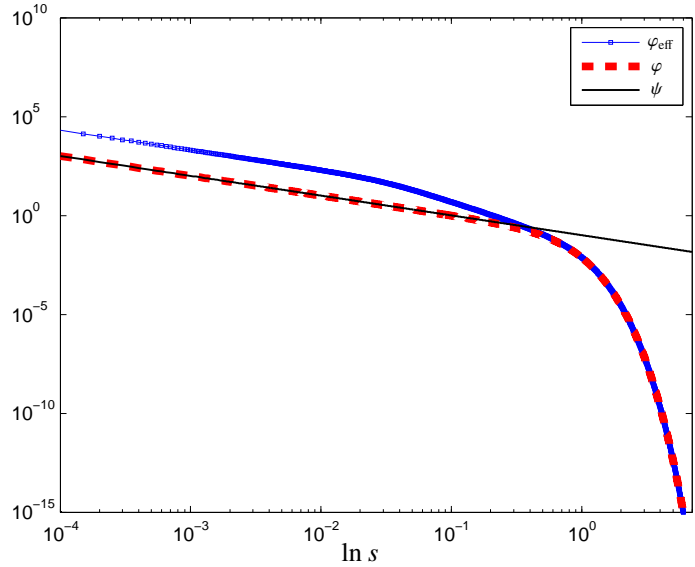

The comparison between the interactions stresses $\varphi_{\textrm{eff}}$, $\varphi$ and $\psi$ is illustrated in Figure \ref{Fig:Interactions}.

From the figure it is clear that $\varphi$ is a good approximation of $\varphi_{\textrm{eff}}$ only for large values of its argument. Therefore the first approximation \eqref{step1} is justified when $h\tilde\rho$ is small. This is exactly the case, though, when the second approximation is not allowed, since $\psi$ and $\varphi$ are close only near zero (corresponding to the opposite case of large $h\tilde\rho$).

Therefore, although this derivation can formally be made, it can not be made rigorous, since the two approximations are mutually incompatible. More precisely, if first truncating to nearest neighbours (replacing $\varphi_{\mathrm{eff}}$ by $\varphi$) and then Taylor-expanding $\varphi$ were a consistent combination, then the Taylor expansion of $\varphi_{\mathrm{eff}}$ and of $\varphi$ would be similar. However, this is not the case, since close to zero $\varphi(t)\approx \psi(t) = 1/\pi^2 t$, while
\footnote{This follows from the two inequalities (since $s\mapsto s\varphi(s)$ is decreasing in $(0,\infty)$ and $\int_{-\infty}^{\infty} s\varphi(s) = \int_{-\infty}^\infty V$)
$$
\varphi_{\mathrm{eff}}(t) = \sum_{k=1}^\infty k \varphi(kt) \leq\sum_{k=1}^\infty \frac{1}{t^2} \int_{(k-1)t}^{kt} s\varphi(s)\, ds =  \frac{1}{t^2}\int_0^\infty s \varphi (s)\, ds,
$$
and
$$
\varphi_{\mathrm{eff}}(t) = \sum_{k=1}^\infty k \varphi(kt) \geq \sum_{k=1}^\infty \frac{1}{t^2} \int_{kt}^{(k+1)t} s\varphi(s)\, ds = \frac{1}{t^2} \int_t^\infty s \varphi(s)\, ds.
$$
}

\begin{equation}
\label{asymp:Veff}
\varphi_{\mathrm{eff}}(t)\approx \dfrac1{2t^2}\int_{-\infty}^{\infty} V \qquad \text{as }t\to0.
\end{equation}
Close to zero, therefore, many neighbours are relevant, and nearest-neighbour truncation results in a large error.

\bigskip

The idea of deriving a continuum dislocation model from a discrete model of dislocation walls in a pile-up is also present in the work of Mesarovic and collaborators \cite{BaskMes10,MesarovicBaskaranetc10}. Their approach consists in performing a two-step upscaling from the discrete dislocation model: first the dislocations are smeared out in the slip plane, but the discreteness is kept inside the wall; then the semi-discrete model obtained in the first step is upscaled in the vertical direction. We stress that there is no physical reason for performing the upscaling in the two directions separately. In fact, as the authors of the quoted papers correctly point out, the error between the predictions of the continuum model they obtain and the discrete model they start with (they refer to it as ``coarsening error") is significant (and this was already observed in \cite{RPG}). Therefore, they need to correct \emph{ad hoc} the resulting continuum model in order to match the predictions of the discrete model.

The continuum model we obtain, instead, agrees perfectly with the discrete model we considered, and therefore we do not have to add any artificial terms to guarantee the matching between the two.

\section{Comments}\label{Comments}
In this section we collect a number of interesting comments on the results presented in the paper.

\medskip

\textit{Transition between the different limit models.}
The transition between the regimes (1)--(5) is \textit{continuous} in terms of the energy.
For example, the integrand in the first term of the continuum energy $E^{(2)}$ in \eqref{E2:c} contains the term $c\,V(cs)$, which for $c$ large (corresponding to a transition from regime (2) to regime (3)) converges to $\left(\int V\right)\delta$, which is indeed the term appearing in the energy $E^{(3)}$ in \eqref{E3:c}. If instead $c$ is small (corresponding to the transition from regime (2) to regime (1)), then the logarithmic singularity of $V$ appears, leading to the energy \eqref{E1:c}.
The other transitions can be explained analogously.

\bigskip

\emph{Mechanical interpretation of $\beta_n$.} The dimensionless parameter $\beta_n$ measures the elastic properties of the medium (described by $K$) in comparison with the strength of the pile-up driving force $\sigma$. Large $\beta_n$, therefore, corresponds to weak forcing, and small $\beta_n$ to strong forcing. Stronger forcing pushes the dislocation walls closer to each other; when $n\beta_n\to 0$, the forcing pushes the walls so close to each other that their distance is always smaller than $h$, so that $V$ is only sampled close to the logarithmic singularity.

Similarly, when $\beta_n\to \infty$, the forcing is so weak that the distance between the walls falls in the exponential tails of $V$.
The intermediate regime is characterized by dislocation walls that span the range from `smaller than $h$' to `larger than $h$'.

\bigskip

\textit{Length scales.} In each of the scaling regimes (1)--(5) it is possible to express the length of the pile-up region in terms of the parameters $K$, $\sigma$, $h$ and the number $n$ of dislocations. As one might intuitively expect, the length of the pile-up region will increase when going from (1) to (5). This is consistent with the previous observation that small $\beta_n$ corresponds to strong forcing and to confining the dislocations in a small domain, while for large $\beta_n$ the forcing is weak, and consequently dislocations can spread out over a larger region.

\smallskip

More precisely, in the scaling regime (1) we found in \eqref{al1} that the length of the pile-up scales like $\frac{K n}{\sigma}$, so in particular it is independent of the wall spacing $h$. This is not surprising, since in this case only the in-plane interactions play a role, so the walls could be equivalently replaced by individual dislocations.
The length of the pile-up in the cases (2), (3) and (4) scales as $\sqrt{\frac{Knh}{\sigma}}$.
The case (5) exhibits an even different length scale, namely $ nh\log \left(\frac{K}{n\sigma h}\right)$.

\medskip

In terms of $\beta_n$ the length of the pile-up region in each scaling regime is given by $\ell_n^{(1)} = n^2\beta_n^2h$, $\ell_n^{(2-4)} = n\beta_n h$ and $\ell_n^{(5)} = \frac{nh}{2\pi}\log\left(\frac{2\beta_n^2}{\pi}\right)$, respectively.

We can explain the transition between the length scales $\ell_n^{(1)}-\ell_n^{(5)}$ in the following way.
We notice that
$$
\ell_n^{(1)} \leq \ell_n^{(2-4)} \quad \Leftrightarrow  \quad \beta_n \leq \frac1n.
$$
Therefore the transition between $\ell_n^{(1)}$ and $\ell_n^{(2-4)}$ happens exactly in the first critical regime (case (2)), where $\beta_n\sim \frac1n$. Concerning the transition between $\ell_n^{(2-4)}$ and $\ell_n^{(5)}$, although $\ell_n^{(5)}\leq \ell_n^{(2-4)}$ always, $\ell_n^{(5)}$ is acceptable only when $\beta_n$ is at least order one (since otherwise $\ell_n^{(5)}\leq 0$, inadmissible). And this condition corresponds to the second critical regime (case (4)).

\bigskip

\textit{More general dislocations arrangements.}
A natural extension of this work is to apply our rigorous upscaling procedure to more general dislocation arrangements. A first direction is the study of \textit{random walls}, namely perfectly straight (possibly finite) dislocation walls where the spacing between dislocations is not constant. As pointed out in
\cite{Hall11} the problem in this more general case is genuinely two-dimensional, and hence much more difficult to treat. We refer also to \cite{ZaiserGroma11} for a related discussion on random walls of dislocations.

\section{Conclusion}\label{Conclusion}

This paper unravels the mechanical response of a system of walls of parallel edge dislocations that move along equidistant parallel slip planes. By implementing a rigorous mathematical limit procedure, we have identified five different parameter regimes, as the number $n$ of walls tends to infinity. These regimes are characterized by the asymptotic behaviour of a single dimensionless parameter.

For each of these regimes we have  identified the limiting internal stress that is generated by the dislocation density, and the structure of a pile-up of walls against a hard obstacle. For two of the regimes the expressions that we obtain were known in the literature, and our results provide new insight by clearly delineating the conditions under which these expressions are valid.
The three other cases are new, and the corresponding behaviour has not been studied before.

Although the analysed wall configuration is highly simplified, the rigorous nature of this work implies that these results can be considered as a benchmark: any (possibly more general) model that describes the behaviour of large numbers of dislocations should reproduce at least the behaviour given by the results of this paper when applied to the corresponding idealised situation.

\bigskip
\bigskip
\centerline{\textsc{Acknowledgments}}
\bigskip
\noindent
The research of L. Scardia was carried out under the project number M22.2.09342 in the framework of the Research Program of the Materials innovation institute (M2i) (www.m2i.nl). The research of M. A. Peletier has received funding from the ITN ``FIRST'' of the Seventh Framework Programme of the European Community (grant agreement number 238702), and from the NWO Vernieuwingsimpuls (VICI scheme). The authors would like to thank Jan Zeman and Andrea Braides for their fruitful comments on the paper. Specifically, it was Jan Zeman who pointed us to the possibility of applying $\Gamma$-convergence to this setup.

\bibliographystyle{model1b-num-names}

\end{document}